\documentclass[10pt]{amsart}
\usepackage[margin=1.25in]{geometry}

\usepackage{amsmath, amscd, amsthm, amssymb}
\usepackage{amsfonts}
\usepackage{mathrsfs}
\usepackage{latexsym, epsfig}

\usepackage[pdftex]{color}
\usepackage{comment}
\usepackage{marginnote}
\usepackage{enumitem}
\usepackage{todonotes}

\usepackage{tikz}
\usepackage{tikz-cd}
\usetikzlibrary{matrix, decorations.pathmorphing, arrows}

\usepackage[colorlinks, citecolor=cyan, linkcolor=magenta]{hyperref}
\usepackage{aliascnt}
\usepackage[nameinlink]{cleveref}

\newcommand{\HH}[0]{\mathbb{H}}
\newcommand{\Z}{\mathbb{Z}}
\newcommand{\ZZ}[0]{\mathbb{Z}}

\newcommand{\R}{\mathbb{R}}

\newcommand{\LL}{{\mathcal L}}
\newcommand{\mc}{\mathcal}

\newcommand{\FF}{\mathcal{F}}

\newcommand{\orb}{\mathcal{O}}

\newcommand{\uc}{\mathscr C}

\newcommand{\stab}{\mathrm{stab}}

\newcommand{\core}{\mathrm{core}}
\newcommand{\fr}{\mathrm{fr}}
\newcommand{\frs}{\mathrm{fr}^s}
\newcommand{\fru}{\mathrm{fr}^u}

\newcommand{\wt}{\widetilde}
\newcommand{\ol}{\overline}
\newcommand{\wh}{\widehat}

\newcommand{\mr}{\mathring}

\renewcommand{\phi}{\varphi}
\renewcommand{\epsilon}{\varepsilon}

\newcommand{\del}{\partial}

\newcommand\tsim{\kern-.4em\sim}

\newcommand\ssm{\smallsetminus}

\newcommand{\mb}{\mathbf}

\newcommand{\define}[1]{\textit{#1}}
\newcommand{\defn}[1]{\define{#1}}

\newtheorem{theorem}{Theorem}[section]

\newaliascnt{lemma}{theorem}
\newtheorem{lemma}[lemma]{Lemma}
\aliascntresetthe{lemma}
\crefalias{lemma}{lemma}

\newaliascnt{proposition}{theorem}
\newtheorem{proposition}[proposition]{Proposition}
\aliascntresetthe{proposition}
\crefalias{proposition}{proposition}

\newaliascnt{corollary}{theorem}
\newtheorem{corollary}[corollary]{Corollary}
\aliascntresetthe{corollary}
\crefalias{corollary}{corollary}

\newaliascnt{conjecture}{theorem}
\newtheorem{conjecture}[conjecture]{Conjecture}
\aliascntresetthe{conjecture}
\crefalias{conjecture}{conjecture}

\newaliascnt{problem}{theorem}

\aliascntresetthe{problem}
\crefalias{problem}{problem}

\newaliascnt{motivation}{theorem}

\aliascntresetthe{motivation}
\crefalias{motivation}{motivation}

\newaliascnt{remark}{theorem}
\newtheorem{remark}[remark]{Remark}
\aliascntresetthe{remark}
\crefalias{remark}{remark}

\newaliascnt{construction}{theorem}
\newtheorem{construction}[construction]{Construction}
\aliascntresetthe{construction}
\crefalias{construction}{construction}

\newtheorem{claim}{Claim}
\newtheorem{definition}{Definition}

\theoremstyle{definition}

\DeclareMathOperator{\intr}{int}

\title[Transverse finiteness conjecture]{The Transverse Finiteness Conjecture 
for pseudo-Anosov flows via universal circles
}
\author{Michael P. Landry and Samuel J. Taylor}
\date{\today}                                           

\begin{document}

\begin{abstract}
We prove Mosher's \emph{Transverse Finiteness Conjecture}, posed in 1996. This states that every taut, finite depth foliation $\mc F$ of an atoroidal $3$--manifold is almost transverse to at most finitely many pseudo-Anosov flows, up to orbit equivalence. Our strategy is to organize the flows that are almost transverse to $\mc F$ using their associated universal circles. This technique also establishes the following rigidity result:
if a finite depth foliation $\mc F$ is almost transverse to a pseudo-Anosov flow $\varphi$ without perfect fits, then $\varphi$ is the unique pseudo-Anosov flow that is almost transverse to $\mc F$, up to orbit equivalence. In fact, in this setting, $\mc F$ has only one minimal universal circle.

The mathematics in this article was developed by the authors without the use of AI.
\end{abstract}
\maketitle

\setcounter{tocdepth}{1}
\tableofcontents


\section{Introduction}
	The theory of codimension one foliations of $3$--manifolds has a long history, going back to seminal work of Reeb, Novikov, and Thurston. The modern treatment prioritizes understanding taut foliations, those whose topology most closely reflects the topology of the ambient $3$--manifold, and the construction and study of such foliations occupies a significant corner of geometric topology. This has been increasingly the case over recent years, especially because of interactions with other important $3$-manifold structures, including pseudo-Anosov flows, hyperbolic geometry, contact structures, and Heegaard Floer homology.
	
\subsection{Motivation}
	Let $M$ be a closed, irreducible, and atoroidal $3$--manifold.
	In this article, we are interested in finite depth foliations of $M$ and their (almost) transverse pseudo-Anosov flows. Such foliations were brought to the forefront by Gabai, who developed powerful tools for producing taut foliations. He proved that for any nontrivial class $\alpha \in H_2(M)$,
	 there is a taut, finite depth foliation of $M$ whose set of compact leaves represents $\alpha$ \cite{Gab83}. 

	Before continuing, let us recall the definition of a finite depth foliation. First, a leaf of a foliation $\mc F$ is depth zero if it is compact. Working inductively, a leaf $\lambda$ of $\mc F$ is depth $k$ if $\ol \lambda \ssm \lambda$ contains leaves of depth at most $k-1$ and at least one leaf of depth equal to $k-1$. The foliation $\mc F$ is itself depth $n$ if all its leaves have depth $\le n$ and it has at least one leaf of depth $n$. If this is the case for some $n\ge 0$, then $\mc F$ is said to be a finite depth foliation.
	
	If $\mc F$ is a foliation of depth zero, i.e. all its leaves are compact, then $\mc F$ is a fibration of $M$ over the circle. Since $M$ is atoroidal, the fibration $\mc F$ packages together a trove of additional structure, largely understood from work of Thurston \cite{thurston1986norm} and Fried \cite{Fri79}. In this case, the monodromy of the fibration has a unique pseudo-Anosov representative, and such a map can be suspended to produce a pseudo-Anosov flow 
$\varphi_{\text{susp}}$
 on $M$ that is transverse to $\mc F$. 
	General pseudo-Anosov flows have also come to occupy their own important corner of geometric topology, especially from the foundational work of Mosher \cite{Mos92,mosher1992dynamical}, Fenley \cite{Fenley2002Rcovered, Fen09,fenley2012ideal}, Calegari \cite{Calegari2000Rcovered, Calegari_promoting}, and more recent breakthroughs of Barthelm\'e, Mann, and their collaborators \cite{BM24, BFM, BBM2025Completing}. 

	When $\mc F$ is a fibration of $M$, $\varphi_{\text{susp}}$ is the \emph{unique} pseudo-Anosov flow transverse to $\mc F$, up to orbit equivalence, and the interactions between $\mc F$ and $\varphi_{\text{susp}}$ have
inspired a huge amount of research, especially in the areas of essential laminations, group actions on the circle, and Cannon--Thurston maps (i.e. equivariant maps from those circles into the $2$-sphere boundary of hyperbolic space). Thus, the questions naturally arise: \emph{do more general (finite depth) foliations admit transverse pseudo-Anosov flows?} And, \emph{how unique are such transverse flows?}
\smallskip	 

The first question was taken up by the celebrated (though largely unpublished\footnote{This is being written and extended by Landry--Tsang \cite{LandryTsangStep1,LandryTsangStep2}.}) work of Gabai--Mosher. They show that for every finite depth foliation $\mc F$ of an atoroidal manifold $M$, there is a pseudo-Anosov flow $\varphi$ that is \emph{almost} transverse to $\mc F$. Here, the \emph{almost} refers to the fact that a minor modification of the flow (called a dynamic blowup) is often needed before the flow can be made transverse. For additional details, see \Cref{sec:background}. Their construction builds an almost transverse flow inductively alongside the finite depth foliation (as in Gabai's original construction) and is quite involved. Importantly, the construction itself has built-in flexibility; that is, at each step choices are made resulting in the possibility of many transverse flows.

In his manuscript on the subject from the 1990s, Mosher posed the following:
\begin{conjecture}[Transverse Finiteness Conjecture; \cite{Mos96}]
Each finite depth foliation on an atoroidal manifold $M$ is almost transverse to at most finitely many pseudo-Anosov flows, up to orbit equivalence.
\end{conjecture}


The goal of this paper is to prove Mosher's conjecture; see \Cref{th:main_tft}. We moreover establish that for the special case of pseudo-Anosov flows \emph{without perfect fits}, introduced and popularized by Fenley, the situation is far more rigid. First, recall that a suspension flow (i.e. the one that is transverse to fibrations over the circle) has the property that no two closed orbits are anti-homotopic. That is, one orbit cannot be homotopic to the reverse of the other. This is because each closed orbit must have positive algebraic intersection with each leaf of the fibration. A general flow is said to have \emph{no anti-homotopic orbits} if this property holds. Fenley proves that a pseudo-Anosov flow has no perfect fits if and only if it has no anti-homotopic orbits (\cite[Theorem B]{Fen16}).

Like the case of suspension pseudo-Anosov flows transverse to a fibration, we show that  if a finite depth foliation $\mc F$ is almost transverse to a pseudo-Anosov flow without perfect fits, then this is the unique pseudo-Anosov flow on $M$ almost transverse to $\mc F$ (see \Cref{th:main_npfs}).

\subsection{Results}
We now more formally summarize our main theorems. The first is to confirm Mosher's Transverse Finiteness Conjecture.
	
\begin{theorem}[Transverse Finiteness Theorem] \label{th:main_tft}
Let $M$ be a closed,  oriented, atoroidal $3$--manifold and let $\mc F$ be a finite depth foliation of $M$. Then there are at most finitely many pseudo-Anosov flows on $M$ that are almost transverse to $\mc F$, up to orbit equivalence. 
%
\end{theorem}

We remark that the well-known Finiteness Conjecture (\cite[Problem 2.31]{K3ProblemList2026}) asserts that any such manifold supports at most finitely many pseudo-Anosov flows, up to orbit equivalence, regardless of whether they are transverse to a foliation. See \Cref{sec:connect} for more on this and the connection to recent work on the subject. The proof of \Cref{th:main_tft}, however, does more than just assert the finiteness of flows almost transverse to $\mc F$. Rather, it organizes these flows according to the \emph{universal circles} they determine for the foliation $\mc F$. For more on this idea, see \Cref{sec:idea}.

\smallskip
The following theorem, however, is distinct from the general Finiteness Conjecture and establishes that pseudo-Anosov flows without perfect fits behave similarly to suspension flows in terms of their rigidity. The theorem also proves \cite[Conjecture 1.4]{landry2026simultaneous} for the class of finite depth foliations.

\begin{theorem}[Transverse rigidity; no perfect fits] \label{th:main_npfs}
Let $M$ be a closed,  oriented, atoroidal $3$--manifold and let $\mc F$ be a finite depth foliation of $M$. 
If $\mc F$ is almost transverse to a pseudo-Anosov flow $\varphi$ without perfect fits, then $\varphi$ is the unique pseudo-Anosov flow on $M$ that is almost transverse to $\mc F$, up to orbit equivalence.
\end{theorem}	

This generalizes the special case of 
depth one foliations transverse to pseudo-Anosov flows, which is mostly due to Junzhi Huang. He proved that if $\mc F$ is a depth one foliation of $M$ and $\phi$ and $\psi$ are pseudo-Anosov without perfect fits, each one transverse to $\mc F$, then $\phi$ and $\psi$ are orbit equivalent \cite{Huang_UC}. Our previous work with Minsky then shows that any pseudo-Anosov flow almost transverse to $\mc F$ must also satisfy the no perfect fits condition \cite{landry2026simultaneous}, thereby implying a slightly weaker version of \Cref{th:main_npfs} for depth one foliations (since Huang's result only applies to transverse, rather than almost transverse flows). 

\smallskip
In fact, the universal circle techniques used to prove \Cref{th:main_npfs} can be modified to establish the following stronger result:

\begin{theorem}[Universal circle rigidity]\label{th:main_ucs}
Let $M$ be a closed,  oriented, atoroidal $3$--manifold and let $\mc F$ be a finite depth foliation of $M$. If $\mc F$ is almost transverse to a pseudo-Anosov flow $\varphi$ without perfect fits, then $\partial \orb_\varphi$ is the unique minimal universal circle of $\mc F$, up to isomorphism.
\end{theorem}

The precise definition of a minimal universal circle, as introduced by Thurston and formalized by Calegari--Dunfield \cite{CalDun_UC}, can be found in \Cref{sec:ucs}. \Cref{th:main_ucs} also provides support for the general relationship between universal circles and transverse flows conjectured in \cite[Conjecture 1.5]{landry2026simultaneous}.

\medskip
Along the way to proving our main theorems, we also establish a number of results that may be of independent interest. These include:
\begin{itemize}
\item \emph{Standard hyperbolic metrics for leaves}. When $\mc F$ is a taut foliation of an atoroidal manifold $M$, the leaves of 
$\wt {\mc F}$ inherit a hyperbolic metric, well-defined up to quasi-isometry. We show that \emph{when $\mc F$ is finite depth}, the corresponding hyperbolic metric on each leaf of $\mc F$ is standard (\Cref{prop:standard}), i.e. it contains no hyperbolic half-spaces. In other words, for any leaf $L$ of $\mc F$, the action $\pi_1(L) \curvearrowright \partial \wt L$ is minimal.

\item \emph{Cooper-Long-Reid for noncompact leaves.} For any leaf $L$ of a foliation $\mc F$ that is transverse to an almost pseudo-Anosov flow $\varphi$, the stable/unstable foliations of $\varphi$ restrict to a pair of transverse foliations $L^{s/u}$ of $L$. We show that as long as $L$ is not a fiber in a fibration of $M$ over the circle,
any ray contained in a leaf of $L^{s/u}$ that stays in a compact set \emph{of $L$} must spiral around a closed leaf of that foliation. The precise statement is given in \Cref{prop:periodic_condition}. This generalizes well-known results of Cooper--Long--Reid \cite{cooper1994bundles} and Fenley \cite{fenley1999surfaces} who prove that when $L$ is a compact leaf of $\mc F$ that is not a fiber, $L^{s/u}$ are \emph{finite foliations}, i.e. every ray spirals around a closed leaf of $L^{s/u}$.

\item \emph{How shadows move as leaves vary.} One main tool in the paper is our previous work with Minsky proving that in the general setting of a foliation $\mc F$ and a transverse almost pseudo-Anosov flow, the boundary of the flow space $\orb_\phi$ of $\phi$ is a universal circle of the foliation $\mc F$.  A key to this picture is relating the ideal geometry of the lifts of leaves of $\wt{ \mc F}$ (the lift of $\mc F$ to the universal cover $\wt M$) to the projections, or \emph{shadows}, of these leaves in $\orb_\phi$.
In \Cref{sec:moving_away}, we establish a general picture of how these shadows evolve as one moves in the leaf space of $\wt{\mc F}$.
\end{itemize}

\subsection{Inductive construction of universal circles} \label{sec:idea}
Although much of the theory we develop applies to arbitrary foliations almost transverse to pseudo-Anosov flows, we ultimately capitalize on the inductive structure of finite depth foliations. The main idea is that each flow $\varphi$ transverse to $\mc F$ determines a universal circle $\partial \orb_\phi$ for $\mc F$, and by work of Barthelm\'e--Frankel--Mann, these circle actions determine the flow $\varphi$, up to orbit equivalence. 

The precise definition of a universal circle is given in \Cref{sec:background}, but here we recall that for a (possibly noncompact) foliated $3$-manifold $X$ with hyperbolic leaves, a universal circle is a single action $\pi_1(X) \curvearrowright \uc$ on a circle that collects together the individual actions $\pi_1(\lambda_X) \curvearrowright \partial \lambda$, where $\lambda$ is a leaf of the foliation lifted to the universal cover $\wt X$ and $\lambda_X$ is its image in $X$. Extending our previous work with Minsky, we show in \Cref{sec:open_saturated} that an almost pseudo-Anosov flow $\varphi$ transverse to $\mc F$ determines a \emph{flow space universal circle} $\uc_X^\phi$ for any open $\mc F$--saturated subspace $X \subset M$. When $X =M$, the universal circle $\uc_X^\phi$ is the flow space boundary $\partial \orb_\phi$.

When $\mc F$ is finite depth, our main theorems are proven by showing that the set of possible flow space universal circles is finite.
Consider the closed subspace $\mc F_{\le k}$ of leaves of depth at most $k$; each component of $M \ssm \mc F_{\le k}$ is an open $\mc F$--saturated set that we call a \define{pocket}. 

The following theorem is proven by induction over the depth of the pocket. 

\begin{theorem}[Finiteness of UCs] \label{th:intro_col}
Suppose that $\mc F$ is a depth $n$ foliation of $M$.
Let $X$ be a pocket of $M \ssm \mc F_{\le n-k}$ for $k =1 , \ldots, n+1$. Then 
there is a finite set 
\[
\uc_1, \ldots, \uc_l
\]
of universal circles of $\mc F|_X$ so that if $\varphi$ is any pseudo-Anosov flow on $M$ almost transverse to $\mc F$, the universal circle $\uc^\varphi_X$ is isomorphic to $\uc_i$ for some $1\le i\le l$. 

Moreover, if there exists a pseudo-Anosov flow without perfect fits almost transverse to $\mc F$, then $l=1$. That is, $\uc^\varphi_X$ does not depend on $\varphi$.
\end{theorem}

The proof is completed in \Cref{sec:collate} and, more than just establishing finiteness, highlights the way in which universal circles can be constructed inductively from the universal circles of lower depth pockets. The base case of the induction, when $k=1$, has the following simple structure: each component of $M \ssm \mc F_{\le n-1}$ is itself a fibration over the circle, and in this case there is a unique (minimal) universal circle $\uc$, which is isomorphic to the ideal boundary of the universal cover of a fiber in the fibration. When $k = n+1$, under the convention that $\mc F_{\le -1}=\emptyset$, \Cref{th:intro_col} implies our main theorems (\Cref{th:main_tft} and \Cref{th:main_npfs}). This uses the fact that, when $M$ is atoroidal, the action $\pi_1(M) \curvearrowright \partial \orb_\phi$ determines $\phi$, up to orbit equivalence \cite[Theorem 1.5]{BFM}. 
Alternatively, this 
can be obtained by combining \cite[Corollary E]{BBM2025Completing}, which states that the action on $\partial \orb_\phi$ determines the action on $\orb_\phi$ and Barbot's well-known result that the action on $\orb_\phi$ determines the flow $\phi$, up to orbit equivalence \cite{barbot1995caracterisation}.

\subsection{Connections to the literature: the Finiteness Conjecture} \label{sec:connect}
There have been significant recent advances made on the general Finiteness Conjecture.
First, there are now many examples of $M$ that admit an arbitrary number of Anosov flows with varying properties; see \cite{BBY17, BM22, CP25b, BY23, BSZ25}.

In the converse direction of bounding the number of (pseudo-)Anosov flows on a fixed manifold, there is recent work \cite{BM24, Zun25, CP25a, barthelme2026uniqueness}. For this paper, most relevant is the result of Barthelm\'e--Tsang--Zung who prove that for a fixed atoroidal $3$-manifold $M$, there are at most finitely many pseudo-Anosov flows without perfect fits on $M$, up to orbit equivalence \cite{BTZ26}. Our \Cref{th:main_npfs} serves to complement their work by showing that at most one of them can be almost transverse to a fixed finite depth foliation. 

\smallskip
The authors have been speaking publicly about the results and methods of this paper since late 2025 \cite{Taylor_TFCtalk_1, Taylor_TFCtalk}. During the final preparation of this paper in Summer 2026, there were two major developments. 
First, Bowden--Colin posted a preprint establishing the Finiteness Conjecture for Anosov flows on atoroidal manifolds as well as for a large class of pseudo-Anosov flows
 \cite{bowen_colin_2026_finiteness}. Second, Gabai and Li posted a proof of the full Finiteness Conjecture, again in the atoroidal setting \cite{gabai2026finiteness}. This clearly implies \Cref{th:main_tft}, although it does not imply \Cref{th:main_npfs} or \Cref{th:main_ucs}. 
 
 We hope that our approach, obtaining the Transverse Finiteness Conjecture as consequence of a structure theory for the universal circles of taut finite depth foliations that also yields these new rigidity results, will be independently interesting.

\subsection*{Acknowledgements}
The authors thank Ellis Buckminster, S\'ergio Fenley, Rachel Roberts, Chi Cheuk Tsang, and Jonathan Zung for helpful conversations.

The completion of this article was supported by the National Science Foundation under Grant No. DMS--2424139, while the authors were in residence at the Simons Laufer Mathematical Sciences Institute in Berkeley, California, during the Spring 2026 semester.
Additionally, Landry was supported by NSF grant DMS--2405453, and Taylor was supported by NSF grant DMS--2503113 and the Simons Foundation.

%
%
%
%
%
%
%
%
%
%
%
%
%
%
%
%

\section{Background}
\label{sec:background}

Here we record the background needed for the paper and refer the reader to other sources for additional details. Throughout, $M$ denotes a closed, irreducible, atoroidal $3$--manifold.

\subsection{Flows, flow spaces, and boundaries}
\begin{figure}
\centering
\includegraphics[height=2in]{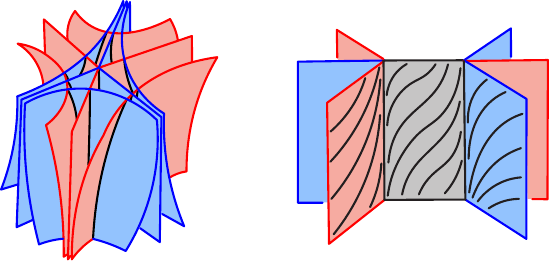}
\caption{Left: a 3-pronged singular orbit of a pseudo-Anosov flow (the flow is upward, the stable foliation is red, and the unstable is blue). Right: a dynamic blowup of a 3-pronged singular orbit whose blowup complex is a single blown annulus (in grey). If more prongs are present, then more complicated blowups are possible.}
\label{fig:blowup}
\end{figure}

For a detailed introduction to the theory of pseudo-Anosov flows, we refer the reader to
the recent book by Barthelm\'e and Mann \cite{barthelme2025pseudo}.
Informally, a \emph{pseudo-Anosov flow} $\phi$ on a closed $3$--manifold $M$ is a flow equipped with a pair of transverse singular $2$-dimensional foliations $W^s$ and $W^u$ in $M$ such that every leaf is a union of flow lines, and every two flow lines on a leaf of $W^s$ are forward asymptotic and every two flow lines on a leaf of $W^u$ are backward asymptotic. There are finitely many singular orbits, each of which has a neighborhood modeled on a semi-branched cover of a closed orbit of an Anosov flow, as in the left side of \Cref{fig:blowup}.

In order to establish our theorems for the natural class of flows transverse to finite depth foliations, we 
must also consider \emph{almost pseudo-Anosov flows}, originally defined by Mosher. An almost pseudo-Anosov flow is the result of replacing finitely many singular orbits with flow-invariant `blowup complexes' in an operation called a \define{dynamic blowup}. See \Cref{fig:blowup}. In this operation, each complex is a connected union of annuli (called \emph{blown annuli}) along their boundary components. Moreover, each blown annulus has the property that its boundary components are closed orbits and every other orbit in the annulus accumulates on distinct boundary orbits in the forward/backward direction.

If $\phi^\sharp$ is a dynamic blowup of a pseudo-Anosov flow $\phi$ obtained in this way, we say $\phi^\sharp$ is an \define{almost pseudo-Anosov flow}. 
There is also a pair of stable and unstable foliations associated to $\phi^\sharp$ that are transverse except along the blowup annuli, which are contained in each foliation (see the right side of \Cref{fig:blowup}). The two flows are related by a blowdown map $M \to M$ that is homotopic to the identity, that collapses each blowup complex to a single singular orbit, and that semiconjugates the flow $\varphi^\sharp$ with its invariant foliations to $\varphi$ with its invariant foliations. For more details see \cite{mosher1992dynamical, Mos96} or our detailed treatment in \cite{LMTstrongtst}. 

We remark that although throughout the paper $\varphi$ will be an almost pseudo-Anosov flow, the reader can first focus on the case where $\varphi$ is an honest pseudo-Anosov flow (sometimes we use the word \emph{honest} to emphasize that a pseudo-Anosov flow really is pseudo-Anosov and not just almost pseudo-Anosov). In most cases, only minor modifications are then needed for the more general case. 

\subsubsection*{The flow space and its boundary}\label{sec:flowspace}

Let $\phi$ be an almost pseudo-Anosov flow on $M$.
Let $\wt \phi$ be the lift of $\phi$ to the universal cover $\wt M$ of $M$.
The stable and unstable foliations $W^s$ and $W^u$ lift to singular foliations $\wt W^s$ and $\wt W^u$ in $\wt M$. A \emph{half-leaf} of a leaf $L$ of $W^{s/u}$ containing the lift of a periodic orbit is the closure of a component of $L-(\text{all lifts of periodic orbits})$.

Let $\orb$ be the quotient of $\wt M$ by the flow lines of $\wt \phi$, which we call the \emph{flow space} of $\phi$ (this is also referred to as the \emph{orbit space} in the literature), and denote the quotient map by $\Theta \colon \wt M \to \orb$.
Fenley--Mosher proved that $\orb$ is a plane \cite[Proposition 4.3]{fenley2001quasigeodesic}. The 2-dimensional foliations $\wt W^s$ and $\wt W^u$ project to 1-dimensional singular foliations of $\orb$ that we denote $\orb^s$ and $\orb^u$, respectively. We call these the stable and unstable foliations of $\orb$. The half-leaves of this foliation are the projections of the half-leaves of $\wt W^{s/u}$.
The action of $\pi_1(M)$ by deck transformations descends to an action of $\pi_1(M)$ on $\orb$ by orientation-preserving homeomorphisms that preserves the stable and unstable foliations.

\begin{figure}
\begin{center}
\includegraphics[height=1.25in]{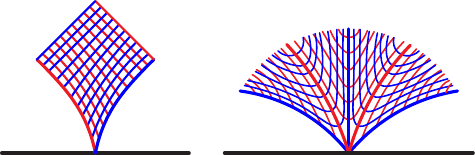}
\caption{In both pictures, the horizontal line denotes $\del\orb$ and $\orb$ lies above it. Left: two foliation rays making a perfect fit, and an associated perfect fit rectangle. Right: any two foliation rays in $\orb$ with the same endpoint in $\del\orb$ are part of a finite sequence of foliation rays ending at that point, each making a perfect fit with the next.}
\label{fig:perfectfits}
\end{center}
\end{figure}

A \emph{perfect fit rectangle} in $\orb$ is the image of a proper embedding of $([0,1]\times[0,1])-(1,1)$ into $\orb$, such that vertical and horizontal lines are mapped into leaves of the stable and unstable foliations. We say the foliation rays corresponding to $\{1\}\times [0,1)$ and $[0,1)\times\{1\}$ \emph{make a perfect fit}, and more generally we say that two foliation rays make a perfect fit if they have subrays lying in the boundary of a perfect fit rectangle as above. See \Cref{fig:perfectfits}. The flow $\varphi$ is said to have \define{no perfect fits} if no foliation rays in $\phi$'s flow space $\orb$ make a perfect fit. As mentioned in the introduction, Fenley shows that $\phi$ has no perfect fits if and only if $\phi$ has no anti-homotopic orbits, and this occurs if and only if each nontrivial $g \in \pi_1(M)$ fixes at most one point of $\orb$ (\cite[Theorem B]{Fen16}). 

In \cite{fenley2012ideal}, Fenley shows that the flow space has a natural compactification
by a circle $\del\orb$ to a closed disk $\ol \orb = \orb \cup \partial \orb$. (Fenley works with non-blown up flows, but his construction is readily adapted to the almost pseudo-Anosov case. Alternatively, one can use the construction in \cite{bonatti2301action}). The action of $\pi_1(M)$ on
$\orb$ extends to an action $\pi_1(M) \curvearrowright \ol \orb = \orb \cup \partial \orb$ by orientation preserving homeomorphisms.
In $\ol \orb$, every foliation ray limits on a unique point, and it is a property of the
construction that the endpoints of foliation rays are dense in $\del\orb$. Moreover, two foliation rays in
$\orb$ have the same endpoint $p\in \del\orb$ if and only if they are part of a finite sequence of rays ending at $p$, each making a perfect fit with the next (\cite[Lemma 3.20]{fenley2012ideal}), see \Cref{fig:perfectfits}.
An important fact about the boundary $\del\orb$ is that if two leaves of $\orb^s$ and/or $\orb^u$ share an ideal point in $\del \orb$, then they are disjoint in $\orb$. An abstract characterization of $\partial \orb$ which works in the almost pseudo-Anosov case and includes all of the required properties can also be found in \cite[Theorem 3]{bonatti2301action}.


\smallskip
We conclude this subsection by discussing a bit more terminology.
A leaf $\ell$ of $\orb^s$ or $\orb^u$ is either homeomorphic to $\R$ or homeomorphic to the interior of a compact connected simplicial tree with at least one vertex of degree $\ge3$ (i.e. the result of deleting the valence 1 vertices from such a compact tree).
In the first case we say that $\ell$ is \emph{regular}, and in the second we say that $\ell$ is \emph{singular}. A point $p\in\orb$ is a \emph{singular point} if it has degree $\ge 3$ in either the stable or unstable leaf containing it
and is otherwise a \emph{regular point}. Alternatively, a regular point is one contained in the interior of a euclidean rectangle such that the restrictions of $\orb^u$ and $\orb^s$ are by horizontal and vertical lines, respectively.

Given a leaf $\ell$ as above, a \emph{leaf face} $l$ in $\ell$ is a properly embedded copy of $\R$ in $\ell$ with the property that at least one of the two components of $\orb -l$ contains no other points of $\ell$. More generally, 
a \emph{leaf slice} of $\ell$ is a properly embedded copy of $\R$ in $\ell$. A point of $\orb$, as well as the leaves containing it, are called \emph{periodic} if it has a nontrivial stabilizer under the $\pi_1(M)$--action. The 
closures of complementary components of the periodic points within a periodic leaf of $\orb^{s/u}$ are called its \emph{half-leaves}. Taking preimages under the quotient $\wt M \to \orb$, we obtain \emph{half-leaves} of $\wt W^{s/u}$, as previously defined.


%

\subsection{Foliations, shadows, and universal circles}
\label{sec:ucs}
Let $\FF$ be a taut foliation of $M$.
That is, $\FF$ is a cooriented foliation of $M$ by surfaces such that there exists a closed curve positively transverse to $\FF$ and intersecting every leaf. By Novikov's theorem \cite{Nov65}, any closed curve transverse to $\FF$ is homotopically nontrivial. Further, each leaf of $\FF$ is $\pi_1$-injective,
 so $\FF$ lifts to a foliation $\wt \FF$ of $\wt M$ by planes.

The quotient of $\wt M$ obtained by collapsing the leaves of $\wt \FF$ is called the \emph{leaf space} of $\FF$, denoted $\LL$. This is a (typically non-Hausdorff) $1$-manifold that also inherits a natural action of $\pi_1(M)$. Typically, we will denote leaves of $\wt {\mc F}$ using greek letters (e.g. $\lambda$, $\mu$, $\nu$) and sometime denote their images in $M$ by $\lambda_M$, $\mu_M$, $\nu_M$. 



\subsubsection*{Transverse foliations and shadows}

We now collect some background on almost pseudo-Anosov flows transverse to foliations, the reference for which is \cite{Fen09}.
We will work in the setting where $\phi$ is an almost pseudo-Anosov flow transverse to a foliation $\FF$ of $M$. In this case, one also says that the corresponding blowdown pseudo-Anosov flow is \emph{almost transverse} to $\FF$. We will always assume that $\phi$ is the \define{minimal} blowup with respect to $\mc F$ in the sense that the foliation on each blowup annulus induced from $\mc F$ is not by segments joining distinct boundary components.

Mosher proved that an (almost) pseudo-Anosov flow on an atoroidal 3-manifold is always transitive, meaning it has a dense orbit \cite[Prop 2.7]{Mos92} (see also \cite[Corollary 1.6]{barthelme2024non}). As such, if $\FF$ is a foliation of $M$ almost transverse to a pseudo-Anosov flow, then it is taut since we can find a closed transversal to all leaves by closing up a long flow segment.

\smallskip
Now let $\wt \FF$ be the lift of $ \FF$ to $\wt M$, and let $\lambda$ be a leaf of $\wt\FF$. The intersection of $\lambda$ with leaves of $\wt W^s$ and $\wt W^u$ defines two singular foliations of $\lambda$ that we denote by $\lambda^s$ and $\lambda^u$.

\smallskip

Let $\Omega_\lambda$ be the projection of $\lambda$ to $\orb$. This open, connected subspace of $\orb$ is called the \emph{shadow} of $\lambda$ and the projection map $\lambda \to \Omega_\lambda$ is a homeomorphism. Each component of the frontier (i.e. topological boundary) $\fr(\Omega_\lambda)$ of $\Omega_\lambda$ in $\orb$ is a stable or unstable leaf slice (see \cite[Theorem 4.1]{Fen09}). 
We denote the union of the stable frontier components by $\frs(\Omega_\lambda)$ and the union of the unstable frontier components by $\fru(\Omega_\lambda)$.
A \emph{frontier chain} of $\Omega_\lambda$ is a sequence of frontier components $\ell_1,\dots, \ell_n$ such that $\ell_i$ shares one ideal point (that is
an ideal point in $\orb \cup \del \orb$) with $\ell_{i+1}$ for $1\le i\le n-1$. The number $n$ is the \emph{length} of the frontier chain. A frontier chain may be comprised of both stable and unstable leaf slices. 
By \cite[Lemma 6.13]{landry2026simultaneous}, the length of all frontier chains is finite (in fact, the length is uniformly bounded above by a constant depending only on $\phi$ and $\mc F$). 

We remark that if the flow $\varphi$ has no perfect fits, then each frontier chain must consist of a single frontier component since in this case no leaf slices share an endpoint in $\partial \orb$.

\begin{figure}
\begin{center}
\includegraphics[height=2in]{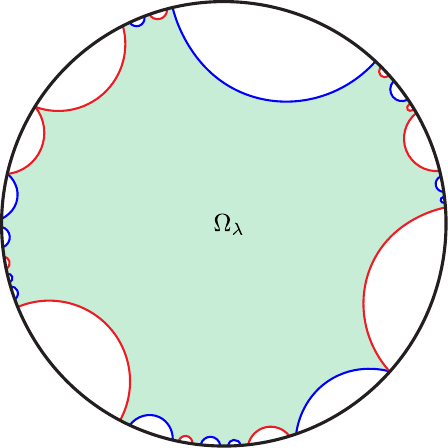}
\caption{The shadow $\Omega_\lambda$ of a leaf $\lambda \subset \wt{\mc F}$.}
\label{fig:shadow}
\end{center}
\end{figure}

The closure of $\Omega_\lambda$ in $\ol\orb$ is a closed disk that we denote by $\ol\Omega_\lambda$ and which we can think of as a union
\[
\ol\Omega_\lambda=\Omega_\lambda\cup \del_\infty \ol\Omega_\lambda\cup \fr(\Omega_\lambda)
\]
where $\del_\infty \ol\Omega_\lambda$ is the \emph{limit set of $\Omega_\lambda$} defined as $\ol \Omega_\lambda \cap \partial \orb$. See \cite{landry2026simultaneous} for details and \Cref{fig:shadow} for an illustration.

\medskip

In \cite[Section 3]{landry2026simultaneous}, we showed that the inverse of the projection homeomorphism $\lambda \to \Omega_\lambda$ extends to a continuous map $\omega_\lambda \colon \ol \Omega_\lambda \to \lambda \cup \partial \lambda$ 
with the properties that
it maps each frontier chain to a single point and its restriction $\partial \ol \Omega_\lambda \to \partial \lambda$ is monotone. 
(This exact formulation is precisely given in \cite[Lemma 5.1]{fenley2026simultaneous}). Here, we recall that a map between circles is \define{monotone} if the preimage of each point is connected.
Since $\omega_\lambda$ is constant on frontier chains of $\Omega_\lambda$, it induces a unique monotone map 
\begin{align}\label{monotone}
\pi_\lambda \colon \partial \orb \to \partial \lambda,
\end{align} 
that is constant on the intervals of $\partial \orb$ that are spanned by frontier chains (see \cite[Construction 3.5]{landry2026simultaneous}). 

\subsubsection*{Universal circles}
\label{sec:ucs}

Let $\mc F$ be a taut foliation of $M$ with leaf space $\LL$.
Since $M$ is {atoroidal}, a theorem of Candel implies that there is a Riemannian metric $g$ on $M$ which pulls back to a hyperbolic metric on each leaf of $\FF$ \cite{Candel_uniformization}; we call such a metric $g$ a \emph{leafwise hyperbolic metric}. 
Hence, each leaf $\lambda$ of $\wt{\mc F}$ is isometric to $\HH^2$  and we can compactify $\lambda$ with its hyperbolic boundary $\partial \lambda$. We define the \defn{circle bundle at infinity} by 
\[
E_\infty=\bigcup_{\lambda\in\mc L}\partial \lambda.
\] 
The bundle $E_\infty$ is topologized so that for any transversal $\tau$ of $\wt {\mc F}$, the natural map 
\[
\text{UT}\wt{\mc F}\mid_\tau\to E_\infty\mid_\tau
\] 
sending a tangent vector $v_p$ to the endpoint of the associated geodesic ray is a homeomorphism.
For additional details, see \cite{CalDun_UC} or \cite{landry2026simultaneous}.

\smallskip
Before defining a universal circle, recall that given a monotone map of circles $m\colon S^1_X\to S^1_Y$, a \emph{gap} of $m$ is a maximal open connected interval in $S^1_X$ sent by $m$ to a single point. The \emph{core} $\mathrm{core}(m)$ of $m$ is the complement of its gaps.

\begin{definition}[Universal circle]\label{def:uc}
	A \define{universal circle} for $\mc F$ is a circle $\uc$ together with an action $\pi_1(M)\curvearrowright \uc$ and a set of monotone maps $\{\pi_\lambda\colon \uc\rightarrow \partial \lambda\:\mid\:\lambda\in\mc L\}$ satisfying:
	\begin{enumerate}
		\item For all $g\in\pi_1(M)$ and for all $\lambda\in\mc L$, the following diagram commutes
		\[\begin{tikzcd}
			{\uc} & {\uc} \\
			{\partial\lambda} & {\partial g\lambda}
			\arrow["{ g}", from=1-1, to=1-2]
			\arrow["{\pi_\lambda}"', from=1-1, to=2-1]
			\arrow["{\pi_{g \lambda}}", from=1-2, to=2-2]
			\arrow["{ g}"', from=2-1, to=2-2]
		\end{tikzcd}\]
		\item The function $\pi\colon \uc \times\mc L\to E_\infty$ defined by \[\pi(z,\lambda)=\pi_\lambda(z)\] is continuous.
		\item For any pair $\lambda$, $\mu$ of incomparable leaves of $\mc L$, the core of $\pi_\lambda$ is contained in a single gap of $\pi_\mu$. \qedhere
	\end{enumerate} \end{definition}

\smallskip	
The universal circle $\uc$ is \defn{minimal} if no open interval in $\uc$ is contained in a gap of $\pi_\lambda$ for all $\lambda \in \mc L$. As in \cite[Lemma 7.1]{landry2026simultaneous}, $\uc$ is minimal if the action $\pi_1(M) \curvearrowright \uc$ is minimal, i.e. every orbit is dense.

\smallskip
Given a universal circle $\uc$, 
if we fix $z\in\uc$, we can define a continuous section $\sigma_z \colon \mc L \to E_\infty$ of $E_\infty$ by the assignment $\sigma_z(\lambda) = \pi_\lambda(z).$ We then set $\mc S_{\uc} = \{\sigma_z : z \in \uc \}$ to be the corresponding set of all such sections of $E_\infty$ determined by $\uc$. 

We say two universal circles $\uc$, $\uc'$ are \emph{isomorphic} if there is a $\pi_1(M)$-equivariant homeomorphism $h\colon \uc\rightarrow\uc'$ such that for any $\lambda \in \LL$, we have $\pi_\lambda^{\uc}=\pi_\lambda^{\uc'}\circ h$. As in \cite[Lemma 3.2]{buckminster2025universal}, two minimal universal circles are isomorphic if and only if their sets of sections of $E_\infty$ are equal. In fact, we will see something even stronger in \Cref{lem:good_order}. Moreover, it follows directly from the definition that if minimal universal circles $\uc$ and $\uc'$ are isomorphic, then the isomorphism $h$ is unique. Indeed, $h$ is uniquely determined on $\core(\pi_\lambda^{\uc})$ for each leaf $\lambda$ and when $\uc $ is minimal, the union of all such cores is dense in $\uc$.

\smallskip


The following theorem of \cite{landry2026simultaneous} is our main source of universal circles in this paper. 

\begin{theorem}[Flow space universal circles]\label{simUC}
Let $M$ be a closed oriented atoroidal manifold and $\phi$ a pseudo-Anosov flow almost transverse to a foliation $\FF$ of $M$.
The action $\pi_1(M) \curvearrowright \del \orb$ and the monotone maps $\pi_\lambda\colon \del\orb\to \del \lambda$ from \cref{monotone}, comprise a universal circle for $\FF$.

Moreover, each gap of $\pi_\lambda \colon \partial \orb \to \partial \lambda$ is spanned by a finite chain of frontier leaves. 
\end{theorem}

\begin{remark}[Universal circles for $\mc F$--saturated sets] \label{rmk:uc_sat}
\emph{Given an open $\mc F$--saturated set $X \subset M$, we can choose a lift $\wt X$ to the universal cover $\wt M$ and consider the circle bundle at infinity restricted to $\wt X$; this determines a circle bundle $E_\infty|_{\wt X} \to \LL|_{\wt X}$. Then a universal circle for the restricted foliation $\mc F|_X$ is precisely as in \Cref{def:uc} with $\pi_1(M)$ replaced by $\pi_1(X)$ and leaves $\lambda$ taken from $\LL|_{\wt X}$. As we will see in \Cref{sec:open_saturated}, an almost transverse pseudo-Anosov flow also naturally determines a universal circle for each open $\mc F$--saturated subspace.}
\end{remark}

\section{Flows and foliations; general structure}
\label{sec:general_structure}

Here we establish some key structural results related to flows and almost transverse foliations. Although our ultimate interest is in finite depth foliations, most of the results in this section do not require the finite depth condition.

\subsection{Spiraling stable/unstable rays}
%
%

Let $\mc F$ be a foliation transverse to an almost pseudo-Anosov flow $\varphi$ that is minimally blown up with respect to $\mc F$. Each leaf $\lambda_M$ of $\mc F$ inherits a pair of singular foliations $\lambda_M^{s/u}$ obtained by intersecting the stable/unstable foliations of $\varphi$ with $\lambda_M$. A \define{stable/unstable ray} in $\lambda_M$ is defined to be an immersed copy of $[0,\infty)$ in a leaf of $\lambda_M^{s/u}$ whose lift to the universal cover is a proper embedding.
Such a ray (viewed as lying in $M$) may be contained in a blown annulus, in which case it is both stable and unstable and we call it a \define{blown ray}.

Each component of intersection between $\lambda_M$ and a blown annulus complex meets at most one periodic orbit of the complex, because $\phi$ is minimally blown up with respect to $\mc F$. Moreover, the component is either
\begin{itemize}
\item a circle, which we call a \define{blown circle}
\item a \define{blown line} containing only interior points of a blown annulus
\item a \define{blown prong} consisting of one or more blown rays with common basepoint, and the basepoint is the intersection of $\lambda_M$ with a boundary orbit of a blown annulus.
\end{itemize}

%

	\begin{proposition}[Spiraling stable/unstable rays]
	\label{prop:periodic_condition}
	Let $\varphi$ be an almost pseudo-Anosov flow transverse to a foliation $\mc F$. Suppose that $\lambda_M$ is a leaf of $\mc F$ that is not a cross section of $\varphi$ and let
 $r$ be a stable/unstable ray in  $\lambda_M$. Then, with respect to the path topology on $\lambda_M$, the ray $r$ is contained in a compact subspace of $\lambda_M$ if and only if it spirals around a closed leaf of the foliation $\lambda_M^{s/u}$.
	\end{proposition}

The proposition is obviously false when $\lambda_M$ is a cross section: in this case $\lambda_M$ is closed, has a well-defined first return map which is pseudo-Anosov, and $\lambda^{s/u}_M$ are the stable/unstable foliations of this pseudo-Anosov. 
In particular, each stable/unstable ray, despite staying in a compact set, is dense in $\lambda_M$. 

\smallskip
For the proof, we use a \emph{veering triangulation} $\tau$  associated to $\varphi$ to represent the stable/unstable ray $r$ combinatorially using the stable/unstable train track determined by $\tau$. The reader is referred to \cite[Section 2]{LMTstrongtst} for the necessary background.

\begin{proof}

Let $r$ be a stable ray in $\lambda_M$. By restricting to a subray, we can assume $r$ contains no singular points. We assume that $r$ does not already parameterize a closed leaf of $\lambda_M^{s/u}$.
	
	We begin by analyzing the situation when $r$ is a blown ray. Then, since $r$ does not meet any singular orbits, there is a blown annulus $A$ in $M$ so that $r$ is contained in the intersection $\lambda_M \cap A$. The 2-dimensional foliation $\mc F$ induces a 1-dimensional foliation of $A$ that we will denote $\mc F_A$. Each leaf of $\mc F_A$ is either a closed curve in $\intr(A)$ or is noncompact and spirals onto such a closed leaf. Since $r$ is not a closed curve by definition,
it must spiral in $A$ onto a closed leaf $c$ of $\mc F_A$. Note that $c$ is the union of all points onto which $r$ accumulates in $M$. However, the local picture in a foliation chart for $\mc F$ around any point in $c$ shows that $r$ does not accumulate on $c$ \emph{in the path topology on $\lambda_M$}, because $r$ does not enter the plaque meeting $c$. We conclude that $r$ accumulates nowhere in $\lambda_M$ and hence escapes compact subsets of $\lambda_M$.

Note that this argument shows that every blown line and blown prong of $\lambda^{s/u}_M$ is properly embedded in $\lambda_M$, a fact we will use going forward.

We now consider the general situation where $r$ has no truncation that is a blown ray.  Since $r$ remains in a compact set of $\lambda_M$ and all blown lines and prongs are properly embedded in $\lambda_M$, $r$ cannot accumulate on any blown line or prong. It may spiral around a blown circle, in which case we are done. \emph{From now on we assume this is not the case, so $r$ does not accumulate on any points in $\lambda_M\cap (\text{blown annuli})$.}


	
 Going forward we assume that $r$ is a stable ray, and the proof in the unstable case is symmetric. The proof is very similar to the one given in \cite[Proposition 8.7]{LMTstrongtst}, which establishes something similar in the case of a compact surface transverse to an almost pseudo-Anosov flow.
	
Let $\kappa$ be a finite collection of periodic orbits together with all singular orbits and all orbits in blowup annuli of $\phi$, and which `kills the perfect fits' of $\phi$. 
(See \cite[Section 3]{LMTstrongtst} for a precise definition.) For our purposes, it suffices to know that this is the condition needed for the existence of a veering triangulation of $M \ssm \kappa$, as employed below.) 
Let $\mr M:=M \ssm \kappa$. The foliation $\mc F$ restricts to a foliation $\mr{\mc F}$ of $\mr M$, and $\mr{\mc F}$ is transverse to $\phi|_{\mr M}$. 
Let $\mr\lambda_M:=\lambda_M\cap \mr M$, which is the result of deleting all singular points, blown circles, blown lines, and blown prongs from $\lambda_M$. 

Note that $r$ is disjoint from $\kappa$ by assumption. We claim that $r$ stays in a compact subset of $\lambda_M\ssm BA$, 
where $BA$ is the subset of $M$ consisting of all blown annuli.
This is because it stays in a compact subset of $\lambda_M$, and does not accumulate on any points in $\lambda_M\cap BA$ by our assumption above. 

There is a veering triangulation $\tau$ of $\mr M=M \ssm\kappa$ associated to $\phi$, and by \cite[Theorem C]{LMT21} we can assume that the 2-skeleton $\tau^{(2)}$, which is naturally a cooriented branched surface, is positively transverse to $\phi$. 
By \cite{Zung_transversefoliations}\footnote{In the cited version of \cite{Zung_transversefoliations}, this is proven only when $\varphi$ is pseudo-Anosov, rather than almost pseudo-Anosov. 
The stronger statement needed here will appear in a forthcoming updated version.}, 
$\mr{\mc F}$ is \emph{carried} by $\tau$ in the sense that it can be split open to a lamination that is carried by $\tau^{(2)}$.
From such a carried position, each leaf inherits an ideal triangulation whose triangles are determined by the faces of $\tau^{(2)}$.

Hence, $\mr\lambda_M$ possesses an ideal triangulation $\triangle(\mr\lambda_M)$, and a train track $\mc V^s(\mr\lambda_M)=\mr\lambda_M\cap B^s$ 
 that, up to isotopy, depend only on its carried position. Here, $B^s$ is the stable branched surface associated to $\tau$; see \cite[Section 4.2]{LMT20}).
The induced triangulation and train track are shown in \Cref{fig:carried}.
For any such carried position, the ray $r$ is carried by $\mc V^s(\mr\lambda_M)$, and visits only finitely many triangles. This is because $r$ stays in a compact subset of $ \lambda_M\ssm BA$, so if it fails to stay in a compact subset of $\lambda_M\ssm \kappa$ it must be because it accumulates in $\lambda_M$ on a finite number of isolated singular points; the corresponding punctures in $\lambda_M-\kappa$ are ideal points of only finitely many triangles.


\begin{figure}[htbp]
\begin{center}
\includegraphics[]{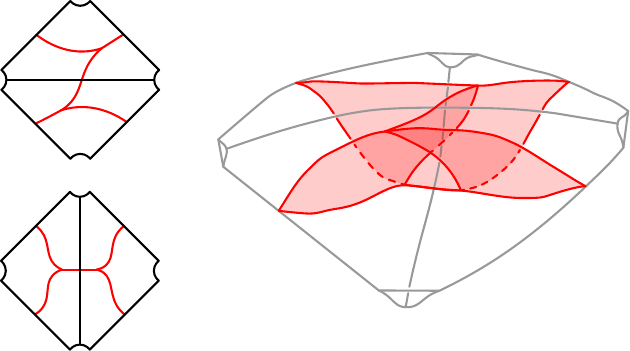}
\caption{On the right, a tetrahedron $t$ of $\tau$ and its intersection with $B^s$. Bottom left: the triangulation and train track on a leaf passing through the bottom two faces of $t$. Top left: the triangulation and train track on a leaf passing through the top two faces of $t$. Each face is foliated by its intersection with $W^{s/u}$.}
\label{fig:carried}
\end{center}
\end{figure}


\emph{Instructive special case.}
When $\lambda_M$ is itself compact, or equivalently $\triangle(\mr\lambda_M)$ has finitely many triangles, the proof uses the same main idea as \cite[Proposition 8.7]{LMTstrongtst}:
essentially, we modify the carried position of $\mr \lambda_M$
by pushing it upward through the tetrahedra of $\tau$ as much as possible, and then analyze the resulting train track. More precisely, whenever $\mr \lambda_M$ traverses the bottom two faces of a tetrahedron, there is a homotopy of $\mr\lambda_M$ along flow lines that sweeps the portion of $\mr\lambda_M$ traversing the bottom faces up through the tetrahedron until it traverses the top two faces. This is called an (upward) \emph{flip move}.
In this situation, we say that the bottom edge of the tetrahedron is \emph{upward flippable}. Given an upward flippable edge $e$, one can perform a flip move as above, and the two triangulations on $\lambda_M$ differ by replacing $e$ in a diagonal exchange. See \Cref{fig:carried}.

Since $\mr\lambda_M$ is not a cross section, hence is not a fiber of $\mr M$,  this process must terminate after finitely many steps by \cite[Proposition 4.5]{Landry2019stable}. This results in a new triangulation of $\mr \lambda_M$ and a new carried position for $\mr\lambda_M$. 
 Note since the homotopies corresponding to flip moves are along flow lines, there is a canonical identification of $\mr \lambda_M$ with its image under the homotopy that respects the stable and unstable foliations. We will take advantage of this identification and continue to use the notation $\mr \lambda_M$ and $r$ to denote the image of these objects under a flip move.

	
In the resulting triangulation of $\mr \lambda_M$, the fact that $\mr \lambda_M$ is in its highest position translates into the fact that the stable track $\mc V^s$ on $\mr \lambda_M$ has no large branches. (A large branch of a train track is a branch for which the switches at both boundary points point inward.) 
Here we are using the combinatorial fact that large branches of $\mc V^s$ are dual to upward flippable edges (see, for example, \cite[Section 4.2]{LMT20}). 
This implies that in this carried position $r$ does not traverse any large branches. In particular, 
the tail of $r$ is contained in a closed `sink' $\sigma$ carried by the train track $\mc V^s(\mr\lambda_M)$, as illustrated in \Cref{fig:sinks}. This implies 
that there are distinct ideal triangles $t_1, \dots, t_n$ 
such that the sequence of triangles traversed by $r$ is eventually $(t_1,\dots, t_n)$ repeating. If we glue $t_i$ to $t_{i+1}$ (mod $n$) along the edge traversed by $r$ and delete an open tubular neighborhood of $\kappa$, then we obtain a closed annulus $A_r$. 

\begin{figure}
\begin{center}
\includegraphics[]{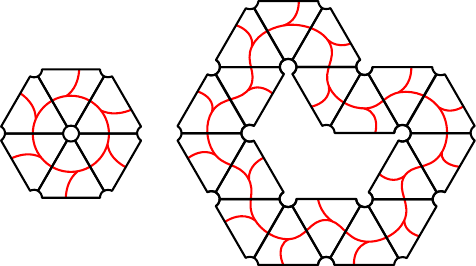}
\caption{The annulus $A_r$ containing the sink $\sigma$.}
\label{fig:sinks}
\end{center}
\end{figure}

We claim that both of the boundary components 
of $A_r$ meet $\mc V^s$, as shown on the righthand side of \Cref{fig:sinks} (rather than in the lefthand side of \Cref{fig:sinks}); equivalently, the sink has branching to both sides. 
Otherwise, the boundary component that is disjoint from $\mc V^s$ corresponds to a puncture of $\lambda_M\ssm \kappa$ that $r$ exits. Since $r$ stays in a compact subset of $\lambda_M\ssm BA$, the puncture corresponds to an isolated singular point of $\lambda^s_M$. However, the train track as in the lefthand side of \Cref{fig:sinks} does not carry the lamination around a singular point, giving a contradiction.
We conclude that there must be a closed leaf of $\lambda^s_M$ onto which $r$ spirals. 


%

\emph{General case.}
In the case that $\lambda_M$ is not compact, the proof is not much different from the special case above. Let $K$ be the union of the finitely many triangles in $\mr\lambda_M$ visited by $r$. 
As in the special case, 
if $K$ contains an upward flippable edge \emph{in its interior}, then we can perform a flip move by pushing $\mr \lambda_M$ up from the bottom to the top of the associated tetrahedron. 
Just as in the special case, each homotopy has the effect of altering $\triangle(\mr \lambda_M)$ by a diagonal exchange on an edge in $\intr(K)$; also, the ray $r$ is still carried by the train track resulting from the flip move.

As in the special case, we repeatedly perform these upward flips along the interior edges of $K$. Since each edge of $\tau$ has finite degree, there can be only finitely many flip moves involving a triangle incident to $\del K$. Applying this reasoning iteratively, we see that the flipping process must terminate after finitely many steps because the initial $K$ is a union of finitely many triangles. The proof finishes as in the special case.
\end{proof}

\subsection{Navigating the leaf space along minimal broken paths}

We will soon detail how shadows vary as we navigate the leaf space $\mc L$. Before thinking about shadows, however, we need to describe {how} to move around in $\mc L$. One does so along ``minimal broken paths," which we now describe.

We use the notation $\lambda \le \nu$ to mean that there is a positive transversal from $\lambda$ to $\nu$ in $\wt M$; this determines a partial order on $\mc L$ by Novikov's theorem. We say two leaves are \define{comparable} if they are related by this partial order, and otherwise they are \define{incomparable}. Note that a transversal between two comparable leaves need not be an orbit of $\varphi$; however, the flow will often be useful as a witness of comparability. For example, for $p \in \orb$, the set of leaves whose shadows contain $p$ 
is totally ordered by $\le$. Alternatively, incomparable leaves of $\wt {\mc F}$ have disjoint shadows in $\orb$.

Since $\wt{\mc F}$ is cooriented and each of its leaves separates $\wt M$, for any leaf $\lambda$ of $\wt {\mc F}$ and any $p \in \wt M \ssm \lambda$, $p$ lies either \define{above} or \define{below} $\lambda$ depending on whether $p$ is on $\lambda$'s positive or negative side. Note that if $\lambda \le \nu$, then $\nu$ is above $\lambda$ and $\lambda$ is below $\nu$. However, the converse is not true.

Two leaves are \define{nonseparated} if they do not have disjoint neighborhoods in the leaf space $\mc L$. A \define{cataclysm} is a collection $\{\nu_i\}$ of incomparable leaves such that there is a monotone sequence $(\lambda_i)_{i \ge0}$ in $\mc L$ such that each $\nu_i$ is accumulated upon by the sequence $(\lambda_i)_{i\ge0}$. In particular, the leaves of a cataclysm are nonseparated in $\mc L$. The cataclysm is from \define{above} or \define{below} depending on whether the $\lambda_i$ are accumulating from above or from below. Said differently, two nonseparated leaves $\nu_1$ and $\nu_2$ are nonseparated from below (i.e. they are part of a cataclysm from below) if and only if each one is below the other.

Because $\mc L$ has a treelike structure, there are unique nonbacktracking `paths' that join arbitrary pairs of leaves in $\mc L$. These are described by Calegari--Dunfield \cite[Section 6.20]{CalDun_UC}, 
where they are called \defn{minimal broken paths}.
More precisely, given $\mu, \nu \in \mc L$, the minimal broken path $\gamma$ between them is the minimal sequence
	\[
	\mu = \nu_0, \nu_1, \nu_2, \ldots, \nu_n = \nu,
	\]
	where $(\nu_{2i}, \nu_{2i+1})$ are ordered (either positively, negatively, or $\nu_{2i}= \nu_{2i+1}$) and 
	$(\nu_{2i+1},\nu_{2i+2})$ are distinct leaves of the same cataclysm. Minimality of the sequence guarantees that the signs of the cataclysms alternate along the path. Hence, after `ordering' the $(\nu_{2i}, \nu_{2i+1})$ for which $\nu_{2i}= \nu_{2i+1}$ appropriately, we see that the ordered intervals also alternate between being positive and negative. We say that the leaves in the positive segments are \defn{ascending} (with respect to $\gamma$) and that the leaves in the negative segments are \defn{descending}. 

	A \define{minimal broken ray} is a semi-infinite sequence of leaves that restricts to a minimal broken path on finite consecutive subsequences.

\subsection{Shadow movement along minimal broken paths} \label{sec:moving_away}
As previously mentioned, a crucial point for us is to understand how the shadows of leaves vary as one moves in the leaf space $\mc L$. Much of our understanding will come from fixing a leaf $\ell$ of $\wt W^{u/s}$ and recording the structure of the induced foliation $\wt{\mc F} \cap \ell$ on $\ell$.

	Recall that for $p\in \orb$, $\gamma_p$ denotes the corresponding orbit in $\wt M$. 
	
	\begin{lemma}[Frontier leaves and order] \label{lem:order}
	Let $p\in \orb$, and let $\lambda$ be a leaf of $\wt \FF$ whose shadow $\Omega_\lambda$ is disjoint from $p$. Let $l$ be the frontier component of $\Omega_\lambda$ that either contains $p$ or that separates $p$ from $\Omega_\lambda$. Then $l$ is a stable frontier component of $\Omega_\lambda$ if and only if $\gamma_p$ lies below $\lambda$. 	The same statement is true replacing ``stable" by ``unstable" and ``below" by ``above."
	\end{lemma}
	
	\begin{proof}
	Suppose that the frontier component $l$ of $\Omega_\lambda$ separating it from $p$ or containing $p$ is stable. Let $l^u$ be a regular unstable leaf meeting $\Omega_\lambda$ and crossing $l$ at some point $q = l \cap l^u$. By considering the intersection of $\ell^u = \Theta^{-1}(l^u)$ with $\lambda$ (see also \cite[Theorem 4.1]{Fen09}), 
we see that $\gamma_q$ (and hence all of $\ell = \Theta^{-1}(l)$) lies below $\lambda$ (i.e. to its negative side). Since $p$ and $q$ lie in the same connected component of $\orb \ssm \Omega_\lambda$, there is a path from $\gamma_p$ to $\gamma_q$ that does not cross $\lambda$. In particular, $\gamma_p$ also lies below $\lambda$. 
	
	A symmetric argument when $l$ is unstable completes the proof.
	\end{proof}

\smallskip

The following lemma will be of fundamental importance later in the paper.

\begin{lemma}[The Moving away lemma]
\label{lem:moving_away}
Let $\lambda \in \mc L$ and let $x \in \partial_\infty \Omega_\lambda \subset \partial \orb$.
Let $\gamma$ be a minimal broken ray in $\mc L$ starting at $\lambda$, and let $\mu$ be a leaf along $\gamma$ such that $x \notin \partial_\infty \Omega_\mu$. Then there is a unique frontier components $l_\mu$ of $\Omega_\mu$ that separates it from $x$. Moreover,

\begin{enumerate}[label=(\alph*)]
\item If $\mu$ is ascending, then $l_\mu$ is stable; if $\mu$ is descending, then $l_\mu$ is unstable.
\item If $\nu$ occurs after $\mu$ along $\gamma$, then $l_\mu$ separates $x$ from $\Omega_\nu$.
\end{enumerate}
\end{lemma}

\begin{proof}
We begin by noting that $\mu$ is ascending (along $\gamma$) if and only if $\lambda$ is below $\mu$ (i.e. to its negative side). Moreover, if $\mu$ is ascending, then all leaves of $\gamma$ past $\mu$ are above $\mu$, and the opposite holds if $\mu$ is descending. These facts are clear by considering the coorientations of leaves along $\gamma$, each of which separates later leaves from $\lambda$. See, for example, \cite[Section 3]{fenley2003pseudo}.

Next, we observe that since $x \in \partial \orb \ssm \partial_\infty \Omega_\mu$, $x$ is separated from $\Omega_\mu$ by some frontier leaf $l_\mu$. Because $x \in \partial_\infty \Omega_\lambda$, which we recall is the set of accumulation points of $\Omega_\lambda$ in $\partial \orb$,
there is some $p \in \Omega_\lambda$ that is also separated from $\Omega_\mu$ by $l_\mu$. Hence, by \Cref{lem:order} and the observations above, $l_\mu$ is stable if and only if $\mu$ is ascending. This proves item $(a)$. 

For item $(b)$, assume that $\mu$ is ascending so that $l_\mu$ is stable. (The other case is similar.) Again by the observations above, if $\nu$ occurs after $\mu$ along $\gamma$, then $\nu$ lies above $\mu$. Since $\Theta^{-1}(l_\mu)$ lies below $\mu$ (as in \Cref{lem:order}), $\mu$ separates $\Theta^{-1}(l_\mu)$ from $\nu$ in $\wt M$. Hence, $\mu$ and $\nu$ lie to the same side of $\Theta^{-1}(l_\mu)$ in $\wt M$. This implies that $\Omega_\mu$ and $\Omega_\nu$ lie to the same side of $l_\mu$ in $\orb$. We conclude that $l_\mu$ also separates $x$ from $\Omega_\nu$ and the proof is complete.
\end{proof}

In essence, the Moving away lemma says that as one moves away from $\lambda$ in the leaf space $\mc L$, the corresponding shadows of leaves in $\orb$ move away from the ideal point $x\in \del \orb$. Crucially, this is true regardless of whether one moves up or down in $\mc L$; however, ascending/descending information can be recovered as described in part (a) of the lemma.

 	\subsection{Structure of foliation saturated subspaces}
	\label{sec:structure}
	
Let $\mr X$ be an open $\mc F$--saturated subset of any foliation $\mc F$ transverse to a flow $\varphi$. Let $X$ be its completion with respect to the path metric, which admits a leaf-preserving immersion back into $M$ (see e.g. \cite[Section 2]{CC_smooth}). 

By Dippolito semistability (see \cite{Dipp} and \cite[Section 2]{CC_smooth}), there is a decomposition
\[
X = Q \cup \left (\Sigma \times I \right),
\]
where $Q$ is a compact sutured manifold, $\Sigma$ is a surface, the restriction of $\mc F$ to $Q$ is a sutured manifold foliation, and orbits of $\phi|_{\Sigma \times I}$ are of the form $x \times I$. This is sometimes called an \defn{octopus decomposition} or \defn{Dippolito decomposition} of $X$.
For additional discussion using the language of sutured manifolds, see  \cite[Section 2]{fenley2001quasigeodesic}.

Note  that the restricted foliation on $\Sigma \times I$ is not necessarily the product foliation. In particular, the foliation on $\partial \Sigma \times I$ is not necessarily a foliation by circles; it could be by lines spiraling on each component of $\partial \Sigma \times \partial I$, or more complicated. The boundary of $X$ consists of collections of leaves $\partial_- X$ and $\partial_+ X$ where $\phi$ points into or out of $X$, respectively. In the same manner, we can define $\partial_\pm Q$, which are compact subsurfaces of the leaves $\partial_\pm X$.

If, in the octopus decomposition, $Q = \emptyset$, then $X$ is called a \define{foliated product}. For any open $\mc F$--saturated subset $U$, at most finitely many components of $U$ are not foliated products (see \cite[Lemma 2.2]{CC_smooth}).


An open $\mc F$--saturated set $\mr X$ (or its completion $X$) is called a \define{mapping torus component} if $\mc F|_{\mr X}$ is a fibration of $\mr X$ over $S^1$. The leaves of such open foliated subsets are called \define{mapping torus leaves}.
We say a mapping torus component $X$ is \define{nontrivial} if it is not a foliated product.

\begin{lemma} \label{lem:nontrivial}
Suppose that the transverse flow $\varphi$ is almost pseudo-Anosov.
If $X$ is a mapping torus component \emph{without} closed orbits, then $X$ is a foliated product and orbits pass from $\partial_- X$ to $\partial_+ X$. Hence, we can rechoose the decomposition of $X$ so that $Q = \emptyset$ and $X = \Sigma \times I$.
\end{lemma}

\begin{proof}
Suppose otherwise. Then there is some point in $\partial_-Q$ whose forward orbit does not meet $\del_+Q$. Since $Q$ is compact, this orbit accumulates on some flow invariant compact set $K$. By the weak closing lemma (e.g. \cite[Proposition 1.4.4]{barthelme2025pseudo}) there is a closed orbit near $K$ and hence contained in $Q$, a contradiction.
\end{proof}

We remark that although the notion of a foliated product \emph{a priori} depends on the particular transverse flow (see \cite{BR_forcing}),
if we restrict to the class of transverse almost pseudo-Anosov flows, then we will see in \Cref{prop:products} that the notion of a foliated product is independent of the flow in the case where $\mc F$ is finite depth.

\subsection{Finite depth foliations have standard leaves}
Here, we restrict our attention to finite depth foliations. Recall that a hyperbolic metric on a surface is \define{standard} if it does not contain an embedded hyperbolic half-plane. Equivalently, the action of the surface's fundamental group on its hyperbolic boundary is minimal.

We will start by proving an elementary result about holonomy of finite depth foliations; we refer the reader to Candel and Conlon's classic textbook \cite[Chapter 2]{candel2000foliations1} for the basic theory of holonomy. For the statement of the lemma, recall that the \defn{$k$th derived subgroup} of a group $G$, denoted $G^{(k)}$, is defined inductively as 
\[
G^{(k)} = [G^{(k-1)},G^{(k-1)}], 
\]
where $G^{(1)} = [G,G]$ is the commutator subgroup.

\begin{lemma}\label{lem:holonomy}
Let $\mc F$ be a foliation of depth $n$, and let $L_0$ be a leaf of $\mc F$ with depth $d_0$. Let $J$ be an arc transverse to $\mc F$ meeting $L_0$ at a single point $x_0$ that is disjoint from any leaf of depth less than $d_0$. If $\gamma$ is a loop in $L_0$ based at $x_0$ and $[\gamma] \in \pi_1(L_0)^{(n-d_0)}$, then there is a subinterval $J_\gamma \subset J$ containing $x_0$ so that for all $x \in J_\gamma$,
the holonomy $h_\gamma(x) \in J$ along $\gamma$ is defined and trivial, i.e.
 $h_\gamma(x) = x$.
\end{lemma}

In particular, the \emph{germinal holonomy} (see \cite[Chapter 2.3]{candel2000foliations1}) is trivial along all elements of $\pi_1(L_0)^{(n-d_0)}$. 

Before giving the proof, we introduce a bit of terminology. Suppose that $L_0$ and $J$ are as in the lemma statement, and that $a$ is a loop in $L_0$ based at $x_0$. We say that $a$ \emph{lifts} to a loop in a leaf $L_1$ based at some $y_0 \in L_1 \cap J$ if there is an annulus $S^1 \times [0,1]$ together with an immersion into $M$ that is transverse to $\mc F$ such that $S^1 \times \{0\}$ maps to $a \subset L_0$, $\{\mathrm{pt.}\} \times [0,1]$ maps to the subarc $J' \subset J$ joining $x_0$ to $y_0$, and $S^1 \times \{1\}$ maps to a `lifted' loop in $L_1$ based at $y_0$. In this case, the lifted loop based at $y_0$ is homotopic along the annulus back to $a$ through loops based at points in $J'$.
  
Note that when $a$ lifts to a leaf $L_1$ based at $y_0$, the holonomy $h_a(y_0)$ is defined and is equal to $y_0$. Further, by considering foliation charts as in the construction of holonomy, it is easy to see that the converse also holds: if the holonomy $h_a(x)$ of a point $x\in J$ along a loop $a$ in $L_0$ is defined and trivial, then $a$ lifts to a (nonunique) loop based at $x$ in the leaf $L$ through $x$. 

\begin{proof}
Let $G = \pi_1(L_0)$. For simplicity, we assume that $[\gamma]$ itself can be written as a $(n-d_0)$-fold iterated commutator using a finite collection of $a_j \in G$. (A general element of the $(n-d_0)$th--derived subgroup will be a word in finitely many such iterated commutators, and our argument is easily generalized.) Realize each $a_j$ as a fixed loop based at $x_0$, and take $J_\gamma \subset J$ sufficiently small so that the holonomy of each point in $J_\gamma$ is defined for each of the finitely many loops formed by taking $k$-fold iterated commutators of the $a_j$ and their inverses for $0 \le k \le n-d_0$. At the expense of again shrinking $J_\gamma$, we also replace $\gamma$ with the (based) homotopic loop formed from the corresponding concatenation of the $a_j$s.

Our goal is to prove the following statement, which implies the lemma after setting $k = n-d_0$:
if $\sigma$ is a loop in $L_0$ based at $x_0$ that is a $k$-fold iterated commutator of the $a_j$s, then for any leaf $L$ of depth at most $d_0+k$ and any $x \in L \cap J_\gamma$, $h_\sigma(x) = x$.

The argument is by induction on $k$.
We first observe that if $L$ is any leaf of $\mc F$ meeting $J_\gamma$ in a finite collection of points, then each $a_j$ lifts to $L$ through any point of $L \cap J_\gamma$. Otherwise, $a_j$ or its inverse would have contracting holonomy from which it would follow that $L$ has a limit point in the interval $J_\gamma$, a contradiction. 
Note that this establishes the base case ($k=0$) since a leaf of depth $d_0$ cannot have infinite intersection with $J$ (this would force a leaf of depth $<d_0$ to intersect $J$, which we are assuming does not happen).

Now assume the statement holds for all nonnegative integers less than $k$, and let $L$ be a leaf of depth $d_0+k$ that intersects $J_\gamma$. 
Fix an arbitrary $x \in L \cap J_\gamma$ and let $\sigma$ be any loop in $L$ that is a $k$-fold iterated commutator in the $a_j$s. If $L \cap J_\gamma$ is finite, then just as above, it follows that $\sigma$ lifts to a loop in $L$ based at $x$. So suppose that this is not the case. Then the set of accumulation points of $L \cap J_\gamma$ is closed and we take $x'$ to be the one closest to $x$ along $J_\gamma$. Note that the leaf $L'$ through $x'$ has depth less than $d_0+k$ and so the interval $[x',x] \subset J_\gamma$ has nonempty interior.

Note that if we write $\sigma = [\sigma_1, \sigma_2]$ where each $\sigma_i$ is a loop in $L_0$ based at $x_0$ that is a $(k-1)$-fold iterated commutator, then the induction hypothesis implies that each $\sigma_i$ has a lift to a loop $\wt \sigma_i$ in $L'$ based at $x'$. We also set $\wt \sigma$ to be the concatenation $[\wt \sigma_1, \wt \sigma_2]$, which itself is a lift of $\sigma$ to $L'$ based at $x'$. In particular, $h_\sigma(x') =x'$. If $h_\sigma(x) \neq x$, then up to replacing $\sigma$ by its inverse, we may assume that the holonomy around $\sigma$ is contracting at $x$, i.e. $h_\sigma(x) \in (x',x)\subset J$. 

By the choice of $x'$, the intersection $L \cap [x',x]$ accumulates only on $x'$. Hence, we can label these points as $x = y_0, y_1, y_2 \ldots$ so that $(y_n)_{g\ge0}$ converges monotonically to $x'$. (This set is infinite because it contains the iterations $h^i_\sigma(x)$ as a subsequence and these converge to $x'$ by assumption.) In this case, the spiraling of $L$ onto $L'$ is determined by a \emph{juncture}, i.e. a homomorphism $\eta \colon \pi_1(L') \to \ZZ$ so that the holonomy of $y_i \in L$ around a loop $\delta \in \pi_1(L')$ is given by $h_\delta(y_i) = y_{i+\eta(\delta)}$, for sufficiently large $i$. 
This is fairly straightforward, but additional details can be found in \cite[Lemma 2.2]{CC_obstruct}. 

In particular, since $\eta$ vanishes on the commutator $\wt \sigma = [\wt \sigma_1, \wt \sigma_2]$ representing an element of $[\pi_1(L'), \pi_1(L')]$, we must have that $h_{\wt \sigma}(y_i) = y_i$ for large $i$. 
Hence, $\wt \sigma$ lifts to a loop in $L$ through $y_i$. But since $\wt \sigma$ was itself a lift of $\sigma$, this implies that $\sigma$ itself lifts to a loop in $L$ through $y_i$.
This implies that $h_\sigma(y_i) = y_i$, contradicting our assumption of contracting holonomy along $\sigma$ in the interval $[x',x) \subset J_\gamma$ and completing the proof.
\qedhere

\end{proof}
	
Using the above lemma we can prove that finite depth foliations in hyperbolic manifolds have standard leaves, which we believe was not previously known.
	
	\begin{proposition}[Standard metrics] \label{prop:standard}
	Let $\mc F$ be a finite depth foliation of an atoroidal manifold $M$, and suppose that $g$ is a leafwise hyperbolic metric. Then the induced hyperbolic metric on each leaf of $\mc F$ is standard. 
	In particular, for each leaf $\mu$ of $\wt{\mc F}$ on the universal cover $\wt M$, the action $\stab(\mu) \curvearrowright \partial \mu$ is minimal. 
	\end{proposition}
		
	\begin{proof}
	Suppose toward a contradiction that the metric on some leaf $L$ of $\mc F$ is nonstandard and let $H$ be the embedded halfspace in $L$ with geodesic boundary. We further choose $L$ so that it has minimal depth among all such leaves; we know that $L$ has depth at least one. Let $H_d \subset H$ be the closed subspace of points at distance at least $d$ (in the hyperbolic metric on $L$) from the boundary of $H$ and let $\ol H_d$ be the closure of $H_d$ in the closed manifold $M$. 
	
	Let $K = \bigcap_{d\ge 0}\ol H_d$. As an intersection of nested compact sets, $K$ is nonempty. Moreover, we claim that $K$ is $\mc F$--saturated; this is similar to the proof that the \emph{limit set} of a leaf is $\mc F$--saturated, as in \cite[Chapter 4]{candel2000foliations1}. This can be seen by considering $x \in K$ and a foliation chart containing $x$. For each $d$, the intersection of $H_d$ with the chart contains plaques (i.e. components of intersection of the chart with leaves of $\mc F$) and these plaques converge to the plaque through $x$. From this the claim follows.
	
	Now let $E$ be a leaf of $K$ of highest depth, which is necessarily less than that of $L$. Since $E$ is a nonelementary hyperbolic surface, there is some loop $\gamma$ in $E$ such that as an element of $\pi_1(E)$ it can be written as a nontrivial $n$--fold commutator, where $n$ is the depth of $\mc F$. Fix $d$ many times larger than the length of $\gamma$ in $E$.
	
	Now fix a transversal $J_\gamma$ as in \Cref{lem:holonomy}, so that for each intersection $x \in L \cap J_\gamma$, $h_\gamma(x) = x$. (Note that such a transversal exists since the set of leaves of depth less than that of $E$ is closed.)
But if $x_i$ is a sequence of such intersections contained in $H_d \cap J_\gamma$ and converging to $x_0 = J_\gamma \cap E$, then $\gamma$ has a lift $\gamma_i$ based at $x_i$. Since the length of $\gamma_i$ is not much more than the length of $\gamma$ (since the metric is continuously varying) and $d$ is much larger than the length of $\gamma$, $\gamma_i$ is contained in the hyperbolic halfspace $H$. This implies that $\gamma_i$ and hence $\gamma$ is nullhomotopic, a contradiction. \qedhere

%
%
%
\end{proof}

The following application will be used in \Cref{sec:collate}. It implies that for a finite depth foliation $\mc F$, the foliated products of $\mc F$ are independent of a chosen transverse almost pseudo-Anosov flow.

\begin{corollary}[Products]\label{prop:products}
Suppose that $\mc F$ is a finite depth foliation. If the completion $X$ of an open $\mc F$--saturated set is a topological product (i.e. $X \cong \partial_- X \times [0,1]$), then for any almost pseudo-Anosov flow $\varphi$ on $M$ transverse to $\mc F$, $X$ is a foliated product. That is, orbits of $\varphi$ pass from $\partial_-X$ to $\partial_+X$.
\end{corollary}

\begin{proof}
We suppose, as we may, that $X$ is connected. Letting $\wt X$ be any lift of $X$ to $\wt M$, there is a homeomorphism $\lambda \times [0,1] \cong \wt X$ taking $\lambda \times \{0\}$ to the leaf $\lambda_0 = \partial_- \wt X$ and $\lambda \times \{1\}$ to the leaf $\lambda_1 = \partial_+ \wt X$. Note that here we are only using that the homeomorphism preserves topological boundary; we are not assuming that it sends the product foliation to $\wt{\mc F}|_{\wt X}$. It suffices to show that all orbits starting at $\lambda_0$ pass through $\lambda_1$, which itself follows from the claim that these leaves have the same shadow in $\orb$; i.e. $\Omega_{\lambda_0} = \Omega_{\lambda_1}$. 

For this, note that the product structure gives that $H: = \stab(\wt X) = \stab(\lambda_0) = \stab(\lambda_1)$. Since $H \curvearrowright \partial \lambda_0$ and $H \curvearrowright \partial \lambda_1$ are minimal by \Cref{prop:standard} (this is the only place where we use that $\mc F$ has finite depth), we have that the $H$--equivariant monotone maps $\pi_{\lambda_0} \colon \partial \orb \to \partial \lambda_0$ and $\pi_{\lambda_1} \colon \partial \orb \to \partial \lambda_1$ must have equal cores. 

Hence, it suffices to know that for any leaf $\lambda$ of $\wt{\mc F}$, its shadow $\Omega_\lambda$ is determined by the core of the associated monotone map $\pi_\lambda$. This follows from the fact that by 
\Cref{simUC}, each gap (a complementary component of the core in $\partial \orb$) of $\pi_\lambda$ is spanned by a finite chain of frontier leaves, each of which is a leaf slice. Moreover, for each pair of endpoints of the gap, there is a unique chain of leaf slices joining them in $\orb$. This is proven for any quasigeodesic, almost pseudo-Anosov flow in \cite[Corollary 6.10]{Fen16}, and any almost pseudo-Anosov flow transverse to a finite depth foliation is quasigeodesic by the main theorem of \cite{fenley2001quasigeodesic}. This completes the proof.
\end{proof}

\section{Open foliation saturated sets, shadows, and universal circles}
\label{sec:open_saturated}
Recall that if $\mc F$ is a finite depth foliation of $M$, then each component of $M \ssm \mc F_{\le k}$ is an $\mc F$--saturated open set and we refer to its completion $X$ as a (depth $k+1$) \define{pocket} of $\mc F$. We also use the same terminology to refer to its lift $\wt X$ in $\wt M$, which is the completion of a component of $\wt M \ssm \wt{\mc F}_{\le k}$.


In the course of our proof, we need to treat these components as `generalized leaves,' and in particular we will consider their shadows in the flow space and associate to them universal circles for the restricted foliation $\mc F|_X$. Ultimately, we will show that this structure can be rebuilt from the foliation itself, up to finitely many choices, but first we must have a detailed understanding of their interactions with transverse flows. Because this is likely to be useful in more general settings, we continue to work with arbitrary foliation-statured open sets of foliations that are transverse to almost pseudo-Anosov flows.
\smallskip


\smallskip
In this section, we fix an almost pseudo-Anosov flow $\varphi$ transverse to an arbitrary foliation $\mc F$. As before, let $\mr X \subset M$ be a connected $\mc F$--saturated open set and let $X$ be its metric completion.
Since the foliation is transversely oriented, the set of boundary leaves of $\partial X$ decomposes into $\partial_- X$ and $\partial_+ X$ as previously described.

There is a natural immersion $X \to M$ (extending the inclusion $\mr X \subset M$) and mapping leaves to leaves. We let $\wt X$ be a lift of $X$ to the universal cover $\wt M$ and
set $\Omega_{\wt X}$ to be the shadow $\Theta(\wt X)$ in $\orb$.

\begin{lemma}[Shadows are shadows] \label{lem:shadows}
 Each frontier component of $\Omega_{\wt X}$ is a frontier component of $\Omega_{\nu}$ for some boundary leaf $\nu \subset \partial \wt X$. Moreover, the stable frontier components of $\Omega_{\wt X}$ are stable frontier components of $\Omega_{\partial_- \wt X}$ and the unstable frontier components of $\Omega_{\wt X}$ are unstable frontier components of $\Omega_{\partial_+ \wt X}$.
 \end{lemma}
 
In the lemma statement,  $\Omega_{\partial_- \wt X}$ is defined as the (necessarily disjoint) union of $\Omega_\nu$ over all boundary leaves $\nu$ of $\partial_-\wt X$, and similarly for $\Omega_{\partial_+ \wt X}$. In fact, by transitivity of the flow, the union $\Omega_{\partial_-\wt X} \bigcup \Omega_{\partial_+\wt X}$ is dense in $\Omega_{\wt X}$. Note that if $X$ is a foliated product, then $\Omega_{\wt X} = \Omega_{\partial_-\wt X} = \Omega_{\partial_+\wt X}$. As in the proof of \Cref{prop:products}, the converse is also true but will not be needed.

 \begin{proof}
 Let $x \in \fr (\Omega_{\wt X})$ and let $\gamma_x$ be the corresponding orbit in $\wt M$. Then $ \gamma_x$ is either above or below $\wt X$, and we assume that it is below $\wt X$ (the case of above is similar). Since each leaf of $\wt {\mc F}$ separates $\wt M$, there is a unique boundary leaf $\nu \subset \partial_- \wt X$ that lies above $\gamma_x$ and separates $\gamma_x$ from $\wt X$. 
 Hence, $x \in \fr (\Omega_\nu)$ and we let $l$ be the frontier component of $\Omega_\nu$ containing $x$. Note that since $\nu$ is above $\gamma_x$, $l$ must be stable (\Cref{lem:order}). It now follows easily that $l$ is the entire frontier component of $\Omega_{\wt X}$ containing $x$ and the lemma follows.
 \end{proof}
 
We note that a stable frontier component of $\Omega_{\wt X}$ might be a frontier component of the shadow of a leaf in $\del_+\wt X$ in addition to $\del_-\wt X$. (And the same is true replacing stable by unstable and $+$ with $-$).


\subsection{Core laminations, universal circles, and monotone maps}
\label{sec:flow_ucs}
For $\mc F$, $\varphi$, and $X$ as above, we now define the \define{core laminations} $\Lambda^\pm_{X, \phi} \subset X$ of $X$. When the transverse flow $\varphi$ is fixed, we also denote these laminations by $\Lambda^\pm_X$.


First, the stable/unstable foliations $W^{s/u}$ of $\varphi$ on $M$ restrict to foliations $W^{s/u}_X$ of $X$. Also, we let $\varphi_X$ denote the restricted semi-flow on $X$, which we consider as a one-dimensional oriented foliation transverse to $\mc F|_X$.
Following \cite[Section 5]{landry2023endperiodic}, the foliations $W^{s/u}_X$ contain singular sublaminations $\Lambda_X^+ \subset W^u_X$ and $\Lambda_X^- \subset W^s_X$ defined as follows. As a subspace of $X$, $\Lambda_X^+$ is the union of all $\varphi_X$--orbits which do not meet $\partial_-X$, and $\Lambda_X^-$ is the union of all $\varphi_X$--orbits which do not meet $\partial_+X$. 
Hence, the components of the intersection $\Lambda_X^+ \cap \Lambda_X^-$ are exactly the orbits of $\varphi$ that are contained in $X$. 

As in \cite[Proposition 5.6]{landry2023endperiodic}, $\Lambda_X^\pm$ are actually singular sublaminations of $W^{u/s}_X$. In more detail, a characterization of the leaves of $\Lambda_X^\pm$ is given by \cite[Lemma 5.4]{landry2023endperiodic}\footnote{In this reference, $X$ is obtained by cutting $M$ along compact leaves, but the same arguments work in the general case.}. If $H$ is a leaf of $W^u_X$, then $H \cap \Lambda_X^+$ is either empty or it is all of $H$, after possibly removing any periodic open half-leaves of $H$ 
that meet $\partial_- X$. Such a half-leaf
only arises as the intersection of $H$ with a blown annulus whose orbits approach the periodic orbit of $H$ in \emph{forward} time.
Moreover, each leaf of $\Lambda_X^+$ is of this form.
A similar leafwise description holds for $\Lambda_X^-$, after replacing $u$ with $s$ and $+$ with $-$.

If $\wt X$ is a lift to the universal cover, we use the notation $\wt \Lambda^\pm_X$ 
for the lifted laminations in $\wt X$. The following lemma characterizes the leaves of these laminations from the perspective of the flow space $\orb$. 

\begin{lemma}[Core laminations in the flow space] \label{lem:core_char}
%

A point $p\in \orb$ is contained in a leaf of the singular lamination $\Theta(\wt \Lambda^+_X) \subset \Omega_{\wt X}$ if and only if it is contained in $\Omega_{\wt X}$ and disjoint from $\Omega_{\del_- \wt X}$. 
Moreover, the endpoint of a ray in a leaf of $\Theta(\wt \Lambda^+_X)$ is contained in $\del \orb$ unless the ray has a truncation contained in a blown segment, in which case its endpoint lies on an unstable component of $\fr(\Omega_{\wt X})$.
\end{lemma}


\begin{proof}
Let $p\in \Omega_{\wt X}$. The statement that $p$ is contained in a leaf of $\Theta(\wt \Lambda^+_X)$ is equivalent to saying that the orbit projecting to $p$ does not cross $\partial_-\wt X$. This is true if and only if $p$ is not contained in $\Omega_{\nu}$ for any boundary leaf $\nu \subset \partial_- \wt X$. 

Next we prove the statement about endpoints. Let $r$ be a ray contained in a leaf of the singular lamination $\Theta(\wt \Lambda^+_X)$. 
Since $\Theta(\wt \Lambda^+_X)$ is a lamination, the endpoint of $r$ lies in $\del \ol{\Omega}_{\wt X}$. Note that the endpoint cannot lie in a stable frontier component, because then $r$ would intersect the shadow of a leaf of $\del_-\wt X$ by \Cref{lem:shadows}. Hence if $r$ is not eventually in a blown segment its endpoint is in $\del \orb$, and if it is eventually in a blown segment then its endpoint must lie in an unstable frontier component of $\fr(\Omega_{\wt X})$.
\end{proof}

We next define a universal circle for the restricted foliation $\mc F|_X$, depending on $\varphi$, that we call the \define{$\varphi$--universal circle} for $X$ 
and denote by
\[
\pi_1(X) \curvearrowright \uc^\varphi_{\wt X}.
\]
For this,  first consider the action $\pi_1(X) \curvearrowright \partial \ol \Omega_{\wt X}$, where $\ol \Omega_{\wt X} \subset \orb \cup \partial \orb$ is a closed disk whose topological boundary is the circle $\partial \ol \Omega_{\wt X}$. By \Cref{lem:shadows}, this circle consists of ideal points in $\del\orb$ as well as frontier components of $\Omega_{\wt X}$ of $\orb^{s/u}$. 
We define 
\[
\uc^\varphi_{\wt X} = \partial {\ol \Omega_{\wt X}} /\sim,
\] 
where $\sim$ collapses each frontier component (and hence each frontier chain) to a point. There is now an induced action of $\pi_1(X)$ on the circle $\uc^\varphi_{\wt X}$, which we define to be the $\phi$-universal circle.

We observe that if $\lambda$ is a leaf of $\wt X$ such that $\Omega_\lambda = \Omega_{\wt X}$, there is a canonical $\pi_1(X)$--equivariant identification $\uc^\varphi_{\wt X} = \partial \lambda$ (see \Cref{sec:ucs}). This is the case, for example, if $\lambda$ is a mapping torus leaf of the component $\wt X$ or if $\lambda$ is any leaf of $\wt X$ and $X$ is a foliated product. Note that if $\wt X'$ is another lift of $X$ to $\wt M$, then the actions of $\pi_1(X)$ on $\uc^\varphi_{\wt X}$ and $\uc^\varphi_{\wt X'}$ are conjugate by a $g \in \pi_1(M)$ mapping $\wt X$ to $\wt X'$.

\smallskip

By construction, there is a naturally induced monotone map
\[
\pi_{\wt X} \colon \partial \orb \to \uc^\varphi_{\wt X} 
\]
obtained by collapsing the arcs that are spanned by frontier components to single points. Again, if $\lambda$ is a leaf of $\wt X$ with $\Omega_\lambda = \Omega_{\wt X}$, then $\pi_{\wt X}  = \pi_\lambda$. 
\smallskip

Finally, the construction of $\uc^\varphi_{\wt X}$  is compatible with the monotone maps $\pi_\nu \colon \partial \orb \to \partial \nu$ in the sense that for any leaf $\nu$ of $\wt X$, possibly a boundary leaf, there is a naturally induced monotone map
\[
\pi^\varphi_{\wt X\nu} \colon \uc^\varphi_{\wt X}  \to \partial \nu,
\]
that factors $\pi_\nu$ through $\pi_{\wt X}$. For this, we have only used 
the fact that $\Omega_\nu \subset \Omega_{\wt X}$ and the definitions of $\pi_{\wt X}$ and $\pi_\nu$.
Since $\nu_M \subset X$ we have that $\stab(\nu) \le \stab(\wt X)$ and so $\pi^\varphi_{\wt X\nu}$  is automatically  $\stab(\nu)$--equivariant. 
In addition this family of maps is $\pi_1(X)$-equivariant, so $\uc^\varphi_{\wt X}$ is a universal circle for the foliation $\mc F|_X$ in the sense of \Cref{rmk:uc_sat}.

In \Cref{prop:minimal}, we will see that, together with these maps, $\uc^\varphi_{\wt X}$ is a minimal universal circle for the restricted foliation $\mc F|_X$. In fact, the action $\pi_1(X) \curvearrowright \uc^\varphi_{\wt X}$ itself is minimal.

Note that we use the notation $\pi^\varphi_{\wt X\nu}$ to emphasize the dependence on $\varphi$.

\begin{remark}[Reconstructing $\uc_{\wt X}$] \label{rmk:frontier_leaves}
\emph{
Our main theorems will follow from an inductive procedure to reconstruct the universal circles $\uc_{\wt X}^\bullet$ from data associated only to $\mc F$, up to finite ambiguity. This process will actually produce the action $\pi_1(X) \curvearrowright \uc_{\wt X}^\bullet$ together with the monotone maps $\pi^\varphi_{\wt X\nu} \colon \uc^\varphi_{\wt X}  \to \partial \nu$ for each leaf $\nu$ of the \emph{interior} of $\wt X$. However, from this universal circle structure it is easy to recover the monotone maps $\pi^\varphi_{\wt X\mu} \colon \uc^\varphi_{\wt X}  \to \partial \mu$ for each boundary leaf $\mu$ of $\partial \wt X$. Indeed, for any sequence $\lambda_i \in \mc L|_{\wt X}$ converging to $\mu$, we have for any $x \in \partial \orb$
\[
\pi_\mu(x) = \lim_{i \ge0} \pi_{\lambda_i}(x),
\]
where the limit takes place in $E_\infty|_{\wt X}$.
}
\end{remark}



Finally, if $\wt Y \subset \wt X$ is a subpocket (i.e. an $\mc F$--saturated open subset), then there is also a naturally defined monotone map
$\uc^\phi_{\wt X} \to \uc^\phi_{\wt Y}$. This also follows from the shadow containment $\Omega_{\wt Y} \subset \Omega_{\wt X}$ exactly as above.

	\subsection{Semi-isolated leaves and the gap finiteness lemma}
	
	We recall from \Cref{sec:structure} that to the open $\mc F$--saturated set $\mr X$, there is an octopus decomposition:
	\[
	X = Q \cup \left (\Sigma \times I \right),
	\]
where $X$ is the metric completion of $\mr X$, $Q$ is a compact sutured submanifold, $\Sigma$ is a surface, 
and orbits on $\Sigma \times I$ are of the form $x \times I$. 

It follows immediately from the definitions that the core laminations are contained in $Q$:
\[
\Lambda^\pm_{X} \subset \mathrm{int}(Q).
\]
This inclusion, and compactness of $Q$, is the key to the finiteness results in this subsection.

First, recall that a leaf of a lamination is \define{semi-isolated} if it is isolated to at least one of its sides (in the leaf space of the lamination).
	
	\begin{proposition}[Semi-isolated leaves are periodic] \label{prop:frontier_periodic}
	Any semi-isolated 
	leaf of the core laminations $\Lambda^\pm_{X}\subset X$ is periodic, and the associated closed orbit is homotopic into $\partial_{\mp} X$ along half-leaves of $\Lambda^\mp_X$. 
	\end{proposition}
	
	The closed leaves $c_p^{u/s}$ of $\Lambda^{\pm} \cap \partial_{\pm} X$ are called the \define{principal unstable/stable loops} of $\partial_{\pm} X$. The proposition and its proof give a bijective correspondence between the principal unstable/stable loops of $\partial_{\pm} X$ and the semi-isolated leaves of the core lamination $\Lambda^{\mp}_X$ (see \Cref{cor:principal_bijection}).
	Note that $c_p^{s/u} \subset \partial_\pm Q$.
	
	Additionally, for any leaf $\mu_M \subset \mc F$, its \define{stable/unstable loops} $c^{s/u}$ are the closed (circular) leaves of $W^{s/u} \cap \mu_M$. When $\mu_M \subset \partial_\pm X$, every principal stable/unstable curve is a stable/unstable loop by definition. However, the converse does not hold if there is an annular leaf of $W^{s/u} \cap X$ such that every orbit of the restriction of $\phi$ to the leaf flows from one boundary component to the other.
	
The following lemma, which we will use to prove \Cref{prop:frontier_periodic} makes explicit how stable and unstable loops are related to the frontier components of shadows.

 \begin{lemma}\label{lem:frontierchainsandloops}
 Let $\mu$ be a leaf of $\wt \FF$, and let $l_1,\dots, l_n$ be a periodic frontier chain of $\Omega_\mu$ fixed by $g\in \stab(\mu)$.  Fix $i$, and suppose that $l_i$ is a stable leaf slice.
\begin{itemize}
\item There is a unique periodic point $p_i\in l_i$ that is the endpoint of an unstable (possibly blown) ray $r$. 
\item Identifying $\mu$ with $\Omega_\mu$, the projection of $r$ to $\mu_M$ is an unstable loop $c_i$.
\item Every other unstable ray in $\Omega_\mu$ terminating at $l_i$ projects to a ray in $\mu_M$ that spirals around $c_i$.
\end{itemize} 
The above is also true switching all instances of stable/unstable.

Furthermore, if $i\ne j$, then $c_i\ne c_j$.
\end{lemma}

For the statement of the lemma it is vital that $g\in \pi_1(\mu_M)$ and not  merely $g\in \pi_1(X)$.

\begin{proof}
In \cite[Sec. 4]{Fen09}, Fenley develops the structure theory of shadows of leaves of foliations transverse to almost pseudo-Anosov flows. His arguments show that a frontier component of a shadow is ``regular to the shadow side." In the case when a frontier component does not contain any blown segments, this amounts to saying that each point in the (stable, say) frontier component is the endpoint of a single unstable ray in the shadow. When the frontier component does contain blown segments, and hence multiple periodic points, there is a unique one that is the endpoint of an unstable ray in the shadow. This was also proved in \cite[Sec 7]{LMTstrongtst} in the special case when the leaf is compact, but the arguments there work in this setting as well (the key point to note in order to apply those arguments here is that because $M$ is compact, every leaf of $\wt \FF$ has a tubular neighborhood such that the restriction of $\wt \phi$ to this neighborhood is a product foliation by intervals).

This implies the first bullet point and, taking the quotient by the group action, the second bullet point follows easily.

Now fix $l_i$, suppose $l_i$ is stable, and let $r_i$ and $c_i$ be as described in the first two bullet points. Let $g_i$ be a group element fixing $p_i$. Given an unstable ray in $\Omega_\mu$ terminating at $l_i$, the $\langle g_i\rangle$-orbit of the ray accumulates on $r_i$. It follows that the image of this ray in $\mu$ spirals around $c_i$. This establishes the third bullet point.

For the final statement, it suffices to note that there is no $h\in \pi_1(M)$ such that $h\cdot l_i=l_j$. By an argument of Fenley (see e.g. \cite[Lemma 6.1]{landry2026simultaneous}) this would imply the existence of a $\Z\oplus\Z$ subgroup of $\pi_1(M)$, contradicting that $M$ is atoroidal.
\end{proof}
	
	\begin{proof}[Proof of \Cref{prop:frontier_periodic}]
	Let $\ell^+$ be a semi-isolated leaf of $\Lambda^+_X \subset Q$. Choose a point $p\in \ell^+$ such that there is a stable half-leaf  containing $p$ and emanating from $\ell^+$ into an isolated side of $\ell^+$. (If $\ell^+$ does not contain any blown half-leaves any $p$ will do; otherwise, $p$ can be chosen to lie on a periodic orbit bounding a blown annulus).
	
	Let $D$ be a small disk transverse to $\phi$ centered at $p$, and let $D^{s/u}$ be the foliations of $D$ induced by $W^{s/u}$. Let $l^u$ be the leaf segment of $D^u$ through $p \in D$ corresponding to the fixed isolated side of $\ell^+$, and let $l^s$ be the half-leaf of $D^s$ through $p$ to the isolated side of $l^u$. 
The backward orbits from points of $l^s \ssm p$ exit $X$ through some component $\mu_M$ of $\partial_-X$. The intersection of the backwards flow saturation of $l^s \ssm p$ with $\mu_M$ is therefore a stable ray $r \subset \mu_M$. We orient $l^s\ssm p$ toward $p$, which induces an orientation on $r$.
Note that since $p$ is in $\ell^+$, its backwards orbit does not meet $\partial_- X$. 
	
	Since $\Lambda^\pm_{X} \subset \mathrm{int}(Q)$, we can take $D$ to lie entirely in $Q$. Since $Q$ is saturated with respect to orbit segments in $X$ by construction, we have that $r \subset \mu_M \cap \partial_- Q$. That is, $r$ is contained in a compact subset of the leaf $\mu_M$. Hence, we can apply \Cref{prop:periodic_condition} to deduce that $r$ spirals in $\mu_M$ around a closed leaf $c^s$ of $\mu_M^s$. Note that if $p$ was already chosen to lie on a periodic orbit (as happens in the blown annulus case), we know in fact that the image of $r$ is equal to $c^s$.
	
	Let $\ell^-$ be the leaf of $W^s_{X} = W^s \cap X$ with $\ell^- \cap \partial_-X = c^s$. We claim that $\ell^-$ is actually a leaf of $\Lambda^-_X$ whose intersection with $\ell^+$ is a closed orbit $\gamma$. Since the portion of $\ell^-$ between $c^s$ and $\gamma$ is an annulus, this will complete the proof.
	
	To prove the claim, we consider shadows in $\orb$. If we choose compatible lifts $\wt X, \wt \ell^+, \wt c^s, \wt r, \mu$ of the objects in $M$, $\Theta(\wt \ell^+)$ is a frontier component of $\Omega_\mu$, and $\Theta(\wt r)$ is a stable ray in $\Omega_\mu$ limiting to a point in $\Theta(\wt \ell^+)$. Note that $\Omega_{\wt X}$ contains both $\Omega_\mu$ and $\Theta(\wt \ell^+)$ by construction.
	The stable loop $c^s$ lifts and projects to a stable leaf segment $\Theta(\wt c^s)$ in $\Omega_\mu$ that is invariant under some $g \in \pi_1(\mu_M)$ and limits to some $g$-stabilized point in the same frontier chain as $\Theta(\wt \ell^+)$ (since $\wt r$ and $\wt c^s$ are asymptotic in $\mu$).
By \Cref{lem:frontierchainsandloops}, $\Theta(\wt c^s)$ terminates on $\Theta(\wt \ell^+)$.




The endpoint of $\Theta(\wt c^s)$ determines a closed orbit $\gamma$ in the flow saturation $\phi^\R(\ell^+)$ of $\ell^+$. The backward orbit of every point in $\phi^\R(\ell^+)$ is eventually contained in $\wt X$, since $\ell^+$ is a leaf of $\Lambda^+_X$. As a consequence $\gamma$ is contained in $X$. Since the forward orbit of every point in $c^s$ is asymptotic to $\gamma$, we see that $c^s$ is contained in a leaf of $\Lambda^-_X$ as desired.
\end{proof}

	\begin{corollary} \label{cor:principal_bijection}
	The leafwise homotopy in \Cref{prop:frontier_periodic} defines a bijection between the principal unstable/stable loops of $\partial_\pm X$ and the isolated sides of the semi-isolated leaves of the core laminations $\Lambda_X^\mp$.
	\end{corollary}
	
\begin{proof}
	The proof of \Cref{prop:frontier_periodic} shows that if $l^+$ is a semi-isolated leaf of $\Lambda^+_X$, then for each isolated face of $l^+$ there is a unique half-leaf $l^-$ of $\Lambda^-_X$ 
 such that $l^-$ contains the closed orbit in $l^+$ and $l^- \cap \partial_-X$ is a principal stable curve. Hence there is a well-defined injective map from the set of isolated faces of $\Lambda^+_X$-leaves to the set of principal stable curves.
	
	We will now show this map is surjective. If $c$ is a principal stable curve of $\partial_-X$, then by definition there is a leaf $l^-$ of $\Lambda^-_X$ so that $l^- \cap \partial_-X = c$. Since $l^- \subset \Lambda^-_X$, forward orbits in $l^-$ are contained in $X$ and since the leaf of $W^s$ containing $l^-$ is periodic (it has nontrivial $\pi_1$), the forward orbits converge to some closed orbit $\omega$. Since all backward orbits in the half-leaf of $l^-$ containing $c$ exit $\partial_-X$, $l^-$ determines a unique isolated face of the leaf $l^+$ of $\Lambda^+_X$ that contains $\omega$. The image of this isolated face (under the map described in the first paragraph) is $c$.
\end{proof}
	
	Given a leaf $\nu$ of $\wt{\mc F}$, a frontier component of $\Omega_\nu$ may or may not be stabilized by a nontrivial element of $\stab(\nu)$. However, the following corollary gives a sufficient condition for this to occur.
	
	\begin{corollary}[Periodic frontier leaves] \label{cor:fontier_periodic}
	Suppose $\nu \subset \partial_- \wt X$ and $l^u$ is a frontier component of $\Omega_\nu$ that is not a frontier component of $\Omega_{\wt X}$. Then $l^u$ is fixed by some nontrivial $g \in \stab(\nu)$. 
	\end{corollary}

	\begin{proof}
	The hypotheses imply that the leaf slice $l^u$ is contained in $\Omega_{\wt X}$. Since shadows of $\Omega_{\partial_- \wt X}$ are disjoint, $l^u$ is not contained in any of these shadows, and so by \Cref{lem:core_char}, $l^u$ is a leaf of $\Theta(\wt \Lambda^+_X)$. In fact, it must be the image of a semi-isolated leaf of $\wt \Lambda^+_X$ because it lies on the frontier of $\Omega_\nu$. Hence, the corollary follows from \Cref{prop:frontier_periodic}.
	\end{proof}

We now can establish the following:

\begin{proposition}[Minimality for $\uc^\phi_{\wt X}$] \label{prop:minimal}
Suppose that the hyperbolic metric on each leaf of $\mc F$ is standard (as is the case when $\mc F$ is finite depth, by \Cref{prop:standard}). Then the action $\pi_1(X) \curvearrowright \uc^{\varphi}_{\wt X}$ is minimal. In particular, $\uc^\varphi_{\wt X}$ is minimal as a universal circle of $\mc F|_{X}$.
\end{proposition}

\begin{proof}
By the standardness assumption, the action $\stab(\nu) \curvearrowright \partial \nu$ is minimal for each leaf $\nu$ of $\wt \FF$. If $\nu$ is contained in $\wt X$, then $\pi_1(X) = \stab(\wt X)$ contains $\stab(\nu)$, and for any  $p \in \uc^\varphi_{\wt X}$, the $\stab(\nu)$--orbit $\stab(\nu)\cdot p$ projects to a dense subset of $\partial \nu$ under $\pi^\varphi_{\wt X\nu}$. Hence, $\ol {\stab(\nu)\cdot p}$ contains the core of $\pi^\varphi_{\wt X\nu}$. Note that this is true for each leaf in $\wt X$, so the closure of $\pi_1(X)\cdot p$ in $\uc^\phi_{\wt X}$ contains the core of each monotone map $\pi^\varphi_{\wt X\nu}$. 

To show that the action $\pi_1(X)\curvearrowright \uc^\phi_{\wt X}$ is minimal, it therefore suffices to show that each interval in $\uc^\phi_{\wt X}$ meets the core of $\pi^\varphi_{\wt X\nu}$ for some leaf $\nu$; in other words, that there is no interval in $\uc^\varphi_{\wt X}$ that is collapsed by each $\pi^\varphi_{\wt X\nu}$.

Suppose for a contradiction that such intervals exist, and let $m_{\wt X}$ be the set of maximal closed intervals in $\uc^\varphi_{\wt X}$ that are collapsed by $\pi^\varphi_{\wt X\nu}$ for all $\wt\FF$-leaves $\nu$ in $\wt X$. This set is $\pi_1(X)$-invariant; hence for any $I\in m_{\wt X}$, and any $g \in \pi_1(X)$, either $g\cdot I=I$ or $gI\cap I=\emptyset$.
Now let 
\[
m_{\orb}=\{\pi_{\wt X}^{-1}(I)\mid I\in m_{\wt X}\}
\]
be the $\pi_1(X)$-invariant collection of intervals of $\partial \orb$ that arise by taking the preimages of the maximal intervals under $\pi_{\wt X} \colon \partial \orb \to \uc_{\wt X}$, and let $[x,y]$ be such an interval. 

Let $[x,y]\in m_\orb$. For each $\wt \FF$-leaf $\nu$ contained in $\wt X$, the set $\partial_\infty \Omega_\nu$ meets $[x,y]$ in at most finitely many points because $\pi_\nu([x,y])$ is a singleton.
On the other hand, $\partial_\infty \Omega_{\wt X} \cap [x,y]$ is uncountable because $\pi_{\wt X}([x,y])$ is a nontrivial arc. 

Now since $\phi$ is transitive, there are points whose forward orbit is dense (see \cite[1.6.9]{fisher2019hyperbolic}) and this implies that $\Omega_{\partial_+ \wt X}$ is dense in $\Omega_{\wt X}$. Hence, we can find an infinite sequence of leaves $\nu_i \subset \partial_+ \wt X$ such that each shadow $\Omega_{\nu_i}$ has a frontier chain $c_i$ that spans an interval $[x_i,y_i] \subset \partial \orb$ containing $[x,y]$ with the property that $x_i \to x$ and $y_i \to y$. Since the shadows of leaves of $\partial_+ \wt X$ are disjoint, these chains do not cross. Hence, the $c_i$ limit to a single chain $c$ of leaves of $\orb^{s/u}$ from $x$ to $y$ in $\ol \orb$. 
However, because the chain $c$ meets $[x,y]$ in a countable set of points, there must be a component of $\Omega_{\wt X} \ssm c$ that is disjoint from the shadows of leaves in $\partial_+ \wt X$. This contradicts the density of $\Omega_{\partial_+ \wt X}$ in $\Omega_{\wt X}$, completing the proof.
\qedhere
\end{proof}

		
\medskip

The following lemma is key to controlling the choices available when reconstructing universal circles from foliation data. 

\begin{lemma}[Gap finiteness] \label{lem:gap_finiteness}
With $\varphi$, $\mc F$, and $X$ fixed as above, let $\mu$ be a boundary leaf of $\wt X$. 

\begin{itemize}
\item If $\mu \subset \partial_\pm \wt X$, then the $\stab(\mu)$--orbits of gaps of $\pi_{\wt X\mu} \colon \uc_{\wt X} \to \partial \mu$ are in natural bijective correspondence with the finite set of homotopy classes 
of \emph{principal} unstable/stable loops in $\mu_M$, which are contained in $\Lambda^{\pm}_X \cap \partial_{\pm} Q$. 

\item For a gap $I$ of $\pi_{\wt X\mu} \colon \uc_{\wt X} \to \partial \mu$, the associated homotopy class of (principal and non-principal) unstable/stable loops is in natural bijection with the stable/unstable frontier components of the frontier chain of $\Omega_\mu$ that spans the gap $\pi_{\wt X}^{-1}(I) \subset \partial \orb$ of $\pi_\mu$. 
\end{itemize}
\end{lemma}

\begin{proof}
To prove the lemma we will describe the underlying bijections.

Fix a nontrivial gap $I$ of $\pi_{\wt X\mu} \colon \uc_{\wt X} \to \partial \mu$, where for concreteness we assume that  $\mu \subset \partial_-\wt X$. Then $\Omega_\mu \subset \Omega_{\wt X}$, where stable frontier leaves of $\Omega_\mu$ are also stable frontier leaves of $\Omega_{\wt X}$ by \Cref{lem:shadows}. Let $J \subset \partial \orb$ be the preimage of $I$ under the monotone map $\pi_{\wt X}$. 
Since $I$ maps to a single point of $\partial \mu$ by hypothesis, $\pi_\mu(J)$ is the same single point of $\partial \mu$, which we call $p$. 

The interval $J$ is spanned by a finite chain of frontier components of $\Omega_\mu$. At least one of these components lies in $\Omega_{\wt X}$ since $\pi_{\wt X}(J) = I$ is not a point, and this component must be unstable by \Cref{lem:shadows}.

We claim that there is a unique primitive $g\ \in \stab(\mu)$ so that $p$ is a $g$--attractor for the action on $\partial \mu$. This is because $J$ is spanned by a finite chain of frontier components of $\Omega_\mu$. Since $\pi_{\wt X}(J) = I$ is not a point, this frontier chain contains some finite collection of unstable leaves $\{l^u_j\}$ that lie in $\Omega_{\wt X}$. (Note that $\{l^u_j\}$ might not contain all the unstable leaves of the frontier chain.)
By \Cref{cor:fontier_periodic} each of these is fixed by some element of $\stab(\mu)$, but since they are joined in a chain, there is a single (primitive) element $g \in \stab(\mu)$ that fixes all of them. By equivariance, this $g$ also fixes $p \in \partial \mu$. 

By \Cref{lem:frontierchainsandloops}, 
 each $l^u_i$ determines a unique stable loop in $\mu_M$ in the homotopy class determined by $g$. 
 By \Cref{cor:principal_bijection}, the \emph{principal} stable loops in this homotopy class are in bijection with $\{l^u_i\}$.
 
This shows that each $\stab(\mu)$-orbit of a gap determines a unique homotopy class of principal unstable/stable curves (depending on whether $\mu \subset \partial_\pm \wt X$). Since all the principal curves are disjoint and contained in the compact surface $\partial_\pm Q = Q \cap \partial_\pm X$, there are only finitely many such homotopy classes. This establishes the first bullet point.

The second bullet point follows directly from \Cref{lem:frontierchainsandloops}. We note that the unstable frontier leaves of $\Omega_\mu$ that span $J$ and which are not among the $l^u_j$ determine non-principal stable loops in the same homotopy class. These will not serve any purpose for us, we mention them simply for completeness.
\end{proof}

\begin{remark}[Positive stabilizers for gaps]
\emph{
As in the above lemma, if $I$ is a nontrivial gap of $\pi_{\wt X\mu} \colon \uc_{\wt X} \to \partial \mu$ whose image in $\partial \mu$ is $p$, we call the primitive $g \in \stab(\mu)$ that has $p$ as an attractor for the action on $\partial \mu$ the \define{positive stabilizer} for $p \in \partial \mu$ or the gap $\pi_{\wt X\mu}^{-1}(p) \subset \uc_{\wt X}$. 
}
\end{remark}

\subsection{Moving away for foliation saturated sets} 
\label{sec:move_again}

The following is an easy consequence of the Moving away lemma (\Cref{lem:moving_away}) in the setting of $\mc F$--saturated open sets. Recall that we have fixed a transverse almost pseudo-Anosov flow $\varphi$.

\begin{lemma}[Moving away; pocket version] \label{lem:moving_saturated}
Let $\lambda$ be a leaf of $\wt {\mc F}$ and fix $x \in \partial_\infty \Omega_\lambda \subset \partial \orb$ that is not the endpoint of any leaf in $\orb^s$ or $\orb^u$. Let $\wt X$ be a pocket of $\wt {\mc F}$ with 
a boundary leaf $\nu \subset \partial \wt X$ that separates $\wt X$ from $\lambda$.  If there is a frontier component $l_x$ of $\Omega_{\nu}$ that separates it from $x$, then $l_x$ is also a frontier component of $\Omega_{\wt X}$.
\end{lemma}

Note that if $\nu \subset \partial_- \wt X$, then the minimal broken path from $\lambda$ to a leaf in $\wt X$ is ascending at $\nu$. Hence, the frontier component $l_x$ in this case is stable (see \Cref{lem:moving_away}).

The following corollary is a reinterpretation of \Cref{lem:moving_saturated} in terms of the associated monotone maps. We remind the reader that the monotone maps depend on our fixed flow $\varphi$, which is absent from the notation for convenience.

\begin{corollary}[The next step] \label{cor:next}
With hypotheses as in \Cref{lem:moving_saturated}, either 
\begin{itemize}
\item $\pi_{\wt X \nu}^{-1}(\pi_\nu(x))$ is a single point, which must be $\pi_{\wt X}(x)$, or
\item there is $g \in \stab(\nu)\subset \stab(\wt X)$ so that $\pi_{\nu}(x) \in \partial \nu$ is the $g$--attractor and $\pi_{\wt X}(x)$ is equal to one of the finitely many $g$--attractors in the gap  $\pi_{\wt X \nu}^{-1}(\pi_\nu(x)) \subset \uc_{\wt X}$.
\end{itemize}
Moreover, if the flow $\varphi$ has no perfect fits, then only the first case occurs.
\end{corollary}

\begin{proof}
Clearly,  $\pi_{\wt X}(x) \in \pi_{\wt X \nu}^{-1}(\pi_\nu(x)) \subset \uc_{\wt X}$, so if the preimage is a single point then the first item holds. So suppose that this preimage is a nontrivial gap (i.e. not a point).

From \Cref{lem:moving_saturated}, we have that $\pi_{\wt X}(x) = \pi_{\wt X}(\partial_\infty l_x)$, where $l_x$ is the frontier component of both $\Omega_\nu$ and $\Omega_{\wt X}$ separating them from $x$. Since the gap is not trivial, there must be another frontier component $l_\nu$ of $\Omega_\nu$ in the same chain as $l_x$ such that $l_\nu \subset \Omega_{\wt X}$. (Note that this cannot occur if $\varphi$ has no perfect fits because then every frontier chain is a single frontier component. This establishes the moreover statement.)
By \Cref{cor:fontier_periodic} there exists an element $g\in\stab(\nu)\le \stab(\wt X)$ fixing $l_\nu$ and hence the entire frontier chain of $\Omega_\nu$ (including $l_x$). Using equivariance of the monotone maps, we see that both $\pi_\nu(x) \in \partial \nu$ and $\pi_{\wt X}(x) \in \uc_{\wt X}$ are fixed by $g$.


To complete the proof, we use the associated Cannon--Thurston map $\mb e\colon \del \orb\to S^2_\infty$. This is a $\pi_1(M)$--equivariant continuous map that identifies the endpoints of each leaf of $\orb^{s/u}$; its existence was shown by 
Fenley \cite{Fen16} and Frankel \cite{frankel2015quasigeodesic}. The reader is referred to \cite{fenley2026simultaneous} for additional details in this setting.

Now, the endpoints of the frontier chain of $\Omega_\nu$ containing $l_x$ all are mapped to a single point $z\in S^2_\infty$ by $\mb e\colon \del \orb\to S^2_\infty$, and since $\mb e$ is $\pi_1(M)$--equivariant, this point is also fixed by $g$. Up to replacing $g$ by $g^{-1}$, we can assume $z$ is an attractor for $g$. 

The equivariance of $\mb e$, along with the fact that $\mb e$ is finite to one (see the proof of \cite[Theorem 1.2]{fenley2026simultaneous}), now implies that each ideal point of $l_x$'s frontier chain is also an attractor for $g$ in $\del \orb$. 
As a consequence, the equivariance of $\pi_{\wt X}$ implies that $\pi_{\wt X}(x) = \pi_{\wt X}(\partial_\infty l_x)$ is an attractor for $g$ in $\uc_{\wt X}$, and the equivariance of $\pi_\nu$ shows that $\pi_\nu(x)$ is the unique $g$--attractor in $\partial \nu$. This completes the proof. \qedhere

\end{proof}

\subsection{Canonical circular orders on sections}
We conclude this section with a discussion of how universal circles are often determined by a single section; see \Cref{cor:section_determines}.

Continuing with the notation from the previous section, let $\mc E_\infty|_{\wt X} \to \mc L|_{\wt X}$ be the restricted circle bundle to the pocket $\wt X$. Since we will be organizing transverse flows using their universal circles, the following observation (implicit in \cite[Section 6.14]{CalDun_UC}) is a key point. 

\begin{lemma}[Unique circular orders]
\label{lem:good_order}
Let $\mc S$ be a set of sections of $\mc E_\infty|_{\wt X} \to \mc L|_{\wt X}$. There is at most one circular order on $\mc S$ that agrees with any universal circle for $\mc F|_X$ whose set of sections contains $\mc S$.

In particular, two nonisomorphic universal circles with minimal actions cannot have a single section in common.
\end{lemma}

\begin{proof}
Let $s_1,s_2,s_3$ be distinct sections of $\mc S$. There are leaves $\lambda$ and $\mu$ of $\wt{\mc F}|_{\wt X}$ such that $s_1(\lambda) \neq s_2(\lambda)$ and $s_2(\mu) \neq s_3(\mu)$, and we let $\gamma \colon I \to \mc L_{\wt X}$ be any \emph{continuous} path 
from $\lambda$ to $\mu$.
Note that if $\lambda$ and $\mu$ are incomparable, this path is necessarily nonembedded, but this will not be an issue.  Let $\gamma^*\mc E_\infty$ be the $\gamma$--pullback of $\mc E_\infty|_{\wt X}$ to $I$, the domain of $\gamma$. This is an oriented circle bundle over the interval, i.e. a cylinder $S^1 \times I$, where $S^1$ is oriented.

Now let $\uc$ be any fixed universal circle for $\mc F|_X$ whose set of sections contain $\mc S$. There is also a circle bundle $\uc \times \mc L|_{\wt X} \to \mc L|_{\wt X}$, from the definition of a universal circle, and we denote the $\gamma$--pullback by $C_\gamma$. Using the continuous fiberwise monotone map $\uc \times \mc L|_{\wt X} \to \mc E_\infty|_{\wt X}$ (from \Cref{sec:ucs}), we get a continuous fiberwise monotone map $C_\gamma \to \gamma^*\mc E_\infty$ between cylinders. The sections $s_1,s_2,s_3$ correspond to three distinct vertical arcs $a_1,a_2,a_3$ in $C_\gamma$ that map to the images of the sections $s_1,s_2,s_3$ in the cylinder $\gamma^*\mc E_\infty$. 

The circular order on $s_1,s_2,s_3$ determined by $\uc$ (note this is independent of $\gamma$) can be recovered in 2 ways: by intersecting $a_1, a_2, a_2$ with a fiber $\uc \times \{pt\}$, or by considering the images of $s_1,s_2,s_3$ in the cylinder $\gamma^*\mc E_\infty$, perturbing them to be disjoint, and taking their order in any fiber circle. The second way does not depend on the universal circle $\uc$, as desired.


The last claim now follows because the orbit of the common section determines an equivariant bijective map between dense subsets of the two universal circles. This map must be circular order preserving by the first part of the lemma, and hence it determines a unique isomorphism.
\end{proof}

Combining \Cref{lem:good_order} and \Cref{prop:minimal} gives:
\begin{corollary} \label{cor:section_determines}
Let $\mc F$ be a finite depth foliation. Suppose that $\phi$ and $\psi$ are almost pseudo-Anosov flows that are transverse to $\mc F$ and let $X$ be a pocket of $\mc F$. If the universal circles $\uc^\phi_{\wt X}$ and $\uc^\psi_{\wt X}$ have a single section of $\mc E_\infty|_{\wt X} \to \mc L|_{\wt X}$ in common, then they are isomorphic. 
\end{corollary}

\section{Mapping torus leaves and monotone maps}\label{sec:top_depth}
Recall that an open $\mc F$--saturated subspace $\mr X \subset M$ or its completion $X$ 
is a \define{mapping torus component} if $\mc F|_{\mr X}$ is a fibration over $S^1$. In this case, we call any leaf $\lambda_M$ of $\mc F|_{\mr X}$ a \define{mapping torus leaf}. We remark that $\mr X$ is the unique component of $M \ssm \left (\ol{\lambda_M} \ssm \lambda_M \right)$ that contains $\lambda_M$. 
We will also refer to the lifts of these leaves to $\wt{\mc F}$ as mapping torus leaves.

As in \Cref{sec:structure}, a mapping torus component $X$ is \define{nontrivial} if it is not a foliated product. In this case, we also say that each of its mapping torus leaves is nontrivial. We apply the same terminology to the lift of $X$ to $\wt M$.

%
%

\smallskip
The connection to transverse flows is the following lemma:

\begin{lemma} \label{lem:return}
A leaf $\lambda_M$ of $\mc F$ is a mapping torus leaf if and only if for any transverse flow $\varphi$ on $M$, there is a first return map $f \colon \lambda_M \to \lambda_M$. 
\end{lemma} 

\begin{proof}
If $f \colon \lambda_M \to \lambda_M$ is a first return map, then the image of the mapping torus in $M$ is an $\mc F$--saturated open set $\mr X$, and clearly it is a mapping torus component.

Conversely, if $\mr X$ is a mapping torus component and a forward orbit through $x \in \lambda_M$ failed to return to $\lambda_M$, then it would accumulate on some minimal set (since $M$ is compact) which is necessarily contained in $\mr X$. Then, just as in Fenley \cite[Lemma 4.2]{fenley1992asymptotic}, a nearly closed trajectory in this minimal set can be perturbed to a closed positive transversal to $\mc F|_{\mr X}$ that doesn't meet $\lambda_M$. 
This contradicts that $\mc F|_{\mr X}$ is a fibration with fiber $\lambda_M$.
\end{proof}

\subsection{Uniquely stabilizing elements}
\label{sec:stab}
Let $X$ be a mapping torus component. Any two mapping torus leaves $\lambda_1,\lambda_2$ in the same lift $\wt X$ of $X$ in $\wt M$ are said to be \define{parallel}. Of course, we can make a similar definition in general, but in this section we are particularly interested in mapping torus leaves. 
Recall that we say that $X$ is nontrivial if it is not a product (\Cref{lem:nontrivial} and \Cref{prop:products}). In this case, we also call the mapping torus leaves associated to $X$ nontrivial.

Given any leaf $\lambda$ of $\wt {\mc F}$, we say that $g \in \pi_1(M)$ \define{uniquely stabilizes} $\lambda$ if $g \lambda = \lambda$ and for any other leaf $\nu$, if $g \nu = \nu$, then $\nu$ is parallel to $\lambda$. In other words, the only leaves of $\wt {\mc F}$ fixed by $g$ are parallel to $\lambda$. 

\begin{lemma}[Unique stabilizers exist]
For any nontrivial mapping torus leaf $\lambda$, there exists $g \in \pi_1(M)$ uniquely stabilizing $\lambda$.
\end{lemma}

Of course, this is false if the mapping torus component $X$ containing $\lambda$ is trivial, i.e. a product, because then every curve in $X$ can be homotoped into $\partial_\pm X$. 

\begin{proof}
Let $X$ be the mapping torus component containing $\lambda_M$ and let $\wt X \subset \wt M$ be its lift that contains $\lambda$.

Let $g\in \stab(\lambda)$. First we claim that if $g$ stabilizes any $\wt \FF$-leaf $\nu$ which is not parallel to $\lambda$, then $g$ fixes some component of $\del \wt X$. Indeed, suppose $\nu$ is an $\wt \FF$-leaf such that $g\nu=\nu$. Then $g$ must fix the unique minimal broken path in $\mc L$ from $\lambda$ to $\nu$, and this path contains a unique component $\mu$ of $\partial_\pm \wt X$. We conclude that $g\mu=\mu$.
 
As a consequence, to prove the lemma it suffices to find $g \in \stab(\lambda)$ that does not fix any component of $\del \wt X$. Equivalently, it suffices to find a loop $l$ in $\lambda_M$ that is not homotopic into a boundary leaf of $X$.


Since $X$ contains a closed orbit (from \Cref{lem:nontrivial}), $\Omega_\lambda \subset \orb$ contains a point $p$ that is fixed under some $h \in \pi_1(X)$. Since $\pi_1(\lambda) \curvearrowright \partial \lambda = \uc_{\wt X}$ is minimal (\Cref{prop:standard}), there is a $g \in \pi_1(\lambda_M) = \stab(\lambda)$ whose fixed points in $\partial \lambda$ are in nonadjacent quadrants determined by the half-leaves at $p$. Here, we are using the standard fact that if the limit set of a Fuchsian group is the circle, then the endpoints of hyperbolic elements are dense in the double boundary.
 It follows that $g$ does not stabilize any component of $\partial_\pm X$ because the shadow of such a component would have to contain $p$, but this is impossible because the corresponding orbit is contained in $X$.
\end{proof}

Now that we know uniquely stabilizing elements exist, the next proposition points to their importance.

\begin{proposition}[Unique stabilizers and $\partial \orb$]\label{prop:north_south}
Let $\varphi$ be an almost pseudo-Anosov
flow transverse to $\mc F$. For any mapping torus leaf $\lambda$ of $\wt{\mc F}$, if $g$ uniquely stabilizes $\lambda$, then $g$ acts with north-south dynamics on $\orb \cup \partial \orb$. In particular, the attractor/repeller of $g$ in $\partial \orb$ is not the endpoint of a leaf.
\end{proposition}

\begin{proof}
We begin with the following claims:
\begin{claim}
Let $\nu$ be a leaf of $\wt{\mc F}$, and suppose $h \in \stab(\nu)$. 
Then $h$ either acts with north-south dynamics on $\orb \cup \partial \orb$ or $h$ translates along a stable/unstable line of $\Omega_{\nu}^{s/u}$.
\end{claim}

Note that the latter condition gives a closed leaf of $\nu_M^{s/u}$.

\begin{proof}
The element $h$ also fixes $\Omega_\nu$. If $h$ were not north-south on $\partial \orb$, then it would fix some $p \in \orb$ (here we are using that $M$ is atoroidal; see \cite[Proposition 3.7]{BFM} or \cite[Proposition 3.1.14]{barthelme2025pseudo}). The orbit corresponding to $p$ does not intersect $\nu$, because $\nu$ is fixed by $h$. Hence $p$ is not in $\Omega_\nu$, so $h$ must fix a frontier component of $\Omega_\nu$ (the one determined by the component of $\ol{\orb-\Omega_\nu}$ containing $p$). The required stable/unstable half-leaf in $\Omega_\nu$ emanates from a periodic point in this frontier component (see \Cref{lem:frontierchainsandloops}).
\end{proof}

\begin{claim}
Let $\lambda$ be a mapping torus leaf of $\wt{\mc F}$, and let $\wt X$ be its lifted mapping torus component. If $g \in \stab(\lambda)$ translates along a leaf line in $\lambda^{s/u}$, then $g$ also stabilizes a boundary leaf of $\del_{\pm} \wt X$.
\end{claim}

\begin{proof}
Suppose $g$ fixes and translates along a stable line in $\lambda^s$. Then $g$ must fix a stable line $l^s$ in $\Omega_{\lambda}^s$, and hence fixes at least one point in $\orb$. Similarly to above, one of these points must lie in $\fr(\Omega_\lambda)$, and in particular on an unstable component $l^u$ of $\fr( \Omega_\lambda)$. By \Cref{lem:shadows}, there is a unique component $\mu$ of $\del_+ \wt X$ such that $l^u\subset \fr( \Omega_\mu)$.
Note that $\Omega_\mu$ meets $l^s$, as does $g \Omega_\mu = \Omega_{g\mu}$. Since distinct components of $\del_+ X$ have disjoint shadows, we must have $\mu = g\mu$.
\end{proof}
Putting the two claims together completes the proof of \Cref{prop:north_south}.
\end{proof}

Finally, we need the following lemma that will allow us to easily plug into the Moving away lemmas (\Cref{lem:moving_away} and \Cref{lem:moving_saturated}).

\begin{lemma}\label{lem:good_to_move}
Suppose $g \in \pi_1(M)$ uniquely stabilizes a mapping torus leaf $\lambda$.  For any transverse almost pseudo-Anosov 
flow $\varphi$ and 
 any leaf $\nu$ of $\wt{\mc F}$ that is not parallel to $\lambda$,
 the shadow $\Omega_\nu$ 
 is separated from each fixed point of $g$ in $\partial \orb$ by its frontier components.
\end{lemma}

\begin{proof}
Let $x = g^{+\infty}_{\orb} \in \partial \orb$. If the lemma were false, then there would be a nonparallel leaf $\nu$ such that $x \in \partial_\infty \Omega_\nu$. Note that $x$ cannot be the endpoint of a frontier component of $\Omega_\nu$ because otherwise $g$ would fix such a leaf and have a fixed point in 
$\orb$, contradicting \Cref{prop:north_south}. 
Hence, $x$ is 
nonisolated in $\del_\infty\Omega_\nu$ from both sides.

Next, note that $\nu$ and $g\nu$ are incomparable leaves. This is because $g$ uniquely stabilizes $\lambda$, so the minimal broken path from $\nu$ to $g\nu$ contains consecutive leaves of $\del_+ \wt X$ (or $\del_-X$) which comprise a cataclysm.


On the other hand, $x$ is a point in both $ \partial_\infty \Omega_\nu$ and $\partial_\infty \Omega_{g\nu}$, where it is nonisolated from \emph{both} sides. This implies that $\Omega_\nu$ and $\Omega_{g\nu}$ intersect nontrivially, a contradiction.
\end{proof}

\subsection{Canonical boundaries and monotone maps $\pi_{\lambda\mu} \colon  \partial \lambda \to \partial \mu$}
Again, fix a mapping torus component $X$ and a mapping torus leaf $\lambda$ of $\wt X$. Since $\mr X$ fibers over $S^1$ with fiber $\lambda_M$, it follows that the action $\pi_1(\lambda_M) \curvearrowright \partial \lambda$ extends canonically to an action $\pi_1(X) \curvearrowright \partial \lambda$. For example, one can use the fact that leaves of the universal cover of $\mr X$ have bounded Hausdorff distance in $\wt M$ to canonically identify their hyperbolic boundaries, from which one can derive the required action. Alternatively, this can be seen using the semidirect product description: $\pi_1(X) = \pi_1(\lambda_M) \rtimes \mathbb{Z}$, where the monodromy action is given by the first return map to $\lambda_M$ as in \Cref{lem:return}.

\smallskip
The following basic lemma gives our initial set of universal circles.

\begin{lemma}[Universal circles of mapping torus components]
\label{lem:MTC}
For any transverse almost pseudo-Anosov flow $\varphi$, the map $\pi_{\wt X\lambda} \colon \uc^\varphi_{\wt X} \to \partial \lambda$ determines a conjugacy of the actions $\pi_1(X) \curvearrowright \uc^\varphi_{\wt X}$ and $\pi_1(X) \curvearrowright \partial \lambda$. 

Moreover, for each leaf $\mu \subset \partial_\pm X$, the monotone map $\pi_{\wt X \mu}\colon \uc^\varphi_{\wt X} \to \partial \mu$ is independent of the flow $\varphi$ (see \Cref{rmk:frontier_leaves}).
\end{lemma}

Henceforth, for $\wt X$ a mapping torus component, we now refer to its universal circle and monotone maps without reference to a flow.
	
%
	
\section{Collating lower depth circles and transverse finiteness}
\label{sec:collate}

The proof of the main theorems proceeds by building a finite collection of universal circles for the finite depth foliation $\mc F$ with the property that the boundary of the flow space for any transverse almost pseudo-Anosov flow is isomorphic to one of these universal circles. 

\subsection{An initial reduction}
To simplify the proof of our main theorems, we first make an initial modification of the finite depth foliation $\mc F$ that does not alter the set of transverse almost pseudo-Anosov flows.  A similar reduction is also done by Fenley--Mosher (\cite[Proposition 2.1]{fenley2001quasigeodesic}) and we refer the reader there for additional details.

First, we recall some terminology. For a finite depth foliation $\mc F$,  let $\mc F_{\le k}$ be the union of leaves of $\mc F$ with depth at most $k$; this is closed subspace of $M$ and hence a sublamination. Then each component of $M \ssm \mc F_{\le k}$ is an $\mc F$--saturated open set and we refer to its completion $X$ as a (depth $k+1$) \define{pocket} of $\mc F$. We also use the same terminology to refer to its lift $\wt X$ in $\wt M$, which itself is defined as the completion of a component of $\wt M \ssm \wt{\mc F}_{\le k}$.

Note that for a depth $k+1$ pocket $X$ the leaves of the restricted foliation $\mc F|_X$ have depth at least $k+1$ in the interior $\mr X$, whereas the boundary $\partial X$ is composed of (finitely many) leaves of depth at most $k$.
Hence, $\mc F_{\le k+1} \ssm \mc F_{\le k}$ determines a lamination on each component of $M \ssm \mc F_{\le k}$ for which each individual leaf is closed (as a subspace of $M \ssm \mc F_{\le k}$) and has depth $k+1$ (as a leaf of $\mc F$).

In general, a depth $k$ leaf $L$ of $\mc F$ is \define{isolated} if it is isolated as a leaf of $\mc F_{\le k}$. 

We say the finite depth foliation $\mc F$ is \define{reduced} if 
\begin{enumerate}
\item each mapping torus pocket of $\mc F$ is nontrivial (see \Cref{sec:structure}), and 
\item  the lamination $\mc F_{\le k+1} \ssm \mc F_{\le k}$ on $M \ssm \mc F_{\le k}$ consists of \emph{finitely} many components, each of which is either a fibration over circle (i.e. a mapping torus component) or a fibration over the interval (i.e. a trivially foliated product).
\end{enumerate}

In what follows, we refer to a trivially foliated product as simply a \define{trivial product}.

	
 
	
\begin{lemma}[Reductions]
\label{lem:reduction}
For any finite depth foliation $\mc F$, we may refoliate product pockets and thicken isolated leaves to obtain a reduced finite depth foliation $\mc F_r$. Each almost pseudo-Anosov flow transverse to $\mc F$ is also transverse to $\mc F_r$, up to orbit equivalence. 
\end{lemma}

\begin{proof}
We begin by recalling the modification made by Fenley--Mosher in the first paragraph of the proof of \cite[Proposition 2.1]{fenley2001quasigeodesic}). The authors there fix an arbitrary flow $\varphi$ transverse to $\mc F$ and begin by collapsing the pockets of $\mc F$ that are foliated products along orbit segments of $\varphi$. The result is a foliation $\mc F'$, still transverse to $\varphi$, having the properties that
\begin{itemize}
\item there are finitely many isolated leaves of $\mc F'$, and
\item if $L$ is a nonisolated leaf of depth $k$, then the component of  $\mc F_{\le k} \ssm \mc F_{\le k-1}$ containing $L$ fibers over the circle with fiber $L$.
\end{itemize}

Moreover, since $\mc F'$ no longer contains pockets which are foliated products, the mapping torus components in the second item are nontrivial.
The authors then thicken each isolated leaf of $\mc F'$ to fibration over a closed interval resulting in a reduced foliation $\mc F''$ that is still transverse to the initial flow $\varphi$.

The only issue is that, in general, the `foliated product' pockets of $\mc F$ as well as the collapsing operation depend on the choice of transverse flow. This is resolved as follows. 
First, choose $\varphi$ to be a fixed almost pseudo-Anosov flow transverse to $\mc F$. Then by \Cref{prop:products}, if $X$ is a foliated product pocket with respect to $\varphi$, then it is a foliated product with respect to \emph{any} almost pseudo-Anosov flow transverse to $\mc F$. Hence, rather than collapse the foliated product pockets, as done above, we can 
\emph{refoliate} $\mc F$ on each foliated product $X$ by replacing $\mc F|_X$ with the product foliation using the identification $X \cong \partial_-X \times [0,1]$ determined by $\varphi$. 

If $\psi$ is a different almost pseudo-Anosov flow transverse to $\mc F$, then the 2-dimensional product foliation induced (via the process above) by $\psi$ on $X$  is isotopic to the one induced by $\phi$, by an isotopy that is constant on $\partial X$ (see \cite[Lemma 6.5]{Wal68}\footnote{Waldhausen's statement is only for $X$ compact, but the argument is easily adapted to the case at hand by using a compact exhaustion.}). 
Said differently, $\psi$ is orbit equivalent via a map isotopic to the identity to a flow which is transverse to the trivially foliated products determined by $\phi$. 
Doing this for each product results in a new foliation $\mc F'''$, still transverse to the same almost pseudo-Anosov flows, up to orbit equivalence, which can be obtained from $\mc F'$ by replacing some of its isolated leaves by trivially foliated products. Hence, if we thicken the remaining isolated leaves of $\mc F'''$ (this does not change the orbit equivalence classes of almost transverse pseudo-Anosov flows) we obtain a foliation $\mc F_r$ which is isotopic to $\mc F''$ and has the required properties.
\end{proof}

Note that if $\mc F$ is reduced, then each pocket $X$ is an \emph{injectively} immersed submanifold of $M$ whose boundary consists of leaves of $\mc F$.

\begin{remark}
\emph{
For the rest of this section, we assume that $\mc F$ is reduced.
}
\end{remark}

\subsection{Inductive structure}
The construction of all flow space universal circles is inductive over the depth of the pocket, starting with top depth components, i.e. components of $M \ssm \mc F_{\le n-1}$. The completion $X$ of each such component is a mapping torus component and \Cref{lem:MTC} establishes the base case of the following inductive claim. For its statement, we set $\mc F_{\le -1} = \emptyset$ so that $M \ssm \mc F_{\le -1} = M$.

As before, for any pocket $X$, let $E_\infty|_{X} \to \mc L|_{X}$ be the restricted ideal circle bundle.

\begin{proposition}[Finiteness of UCs] \label{prop:first_col}
Let $X$ be a pocket of $M\ssm \mc F_{\le n-k}$ where $1\le k\le n+1$. Then 
there is a finite set 
\[
\uc_1, \ldots, \uc_l
\]
of universal circles of $E_\infty|_{X} \to \mc L|_{X}$ so that if $\varphi$ is any pseudo-Anosov flow on $M$ almost transverse to $\mc F$, the universal circle $\uc^\varphi_X$ is isomorphic to $\uc_i$ for some $1\le i\le l$. 

Moreover, if there exists an almost transverse pseudo-Anosov flow without perfect fits, then $l=1$. That is, $\uc^\varphi_X$ does not depend on $\varphi$.
\end{proposition}

The proof is by induction on $k$. As previously observed, the base case where $k=1$ is a consequence of \Cref{lem:MTC}. At the end of this section we show how the case $k=n+1$ proves the main theorems of the paper.

\subsection{Setting up the inductive hypothesis; the dual tree}
\label{sec:setup}
Now suppose that \Cref{prop:first_col} holds up to $1\le k< n+1$; we aim to establish the statement for $k+1$. Toward this goal, let $\mr X$ be a component of $M \ssm \mc F_{\le n-(k+1)}$, and let $X$ be its completion. If $X$ is a mapping torus component, then we are again done by \Cref{lem:MTC}, and so we suppose that this is not the case.
We will apply the inductive hypothesis to each of the finitely many components of $\mr X \ssm \mc F_{\le n-k}$, but first we need an effective way of organizing how these components fit together across depth $n-k$ leaves.

%
%
%
%
%
%

\smallskip

Returning to the case at hand, in the open set $\mr X$, the lamination $\mr X \cap \mc F_{\le n-k}$ (consisting of depth $n-k$ leaves) is a finite collection of trivial products, each of which is embedded in $\mr X$ as a closed subspace. (Here, we are using that $\mc F$ is reduced and that $\mr X$ is not a mapping torus).
Hence, we can construct its directed dual graph $G_X$ in the usual sense: the vertices of $G_X$ are components of $\mr X \ssm \mc F_{\le n-k}$, with a directed edge from $v$ to $w$ if there is a trivial product of
depth $n-k$ leaves whose bottom leaf belongs to $\partial_+ v$ and whose top leaf belongs to $\partial_- w$. Here, 
$\partial_\pm v$
 denotes the positive/negative boundary leaves of the completion $\ol v$ of the component $v$, and bottom/top is with respect to the coorientation on leaves.
Since $\mc F$ is reduced, $G_X$ is a finite graph that we can consider as naturally embedded in $X$; we denote its preimage in $\wt X$ by $T_X$, which we call the \define{dual tree} of depth $n-k$ leaves.

The claim that $T_X \subset \wt X$ is indeed a (simplicial) tree can be verified directly. First, note that there is an induced map $T_X \to \mc L|_X$ from which it easily follows that $T_X$ has no embedded cycles. Further, $T_X$ is connected because for any $\wt v, \wt w$ in $T_X^0$ a continuous path between them in $\wt X \ssm \partial_\pm \wt X$ crosses at most finitely many trivial products of depth $(n-k)$ leaves since such products cannot accumulate \emph{in} $\wt X \ssm \partial_\pm \wt X$. 

\begin{remark}[Notation; indexing pockets using $T_X$]
\emph{
In what follows, we often denote vertices of $T_X$ by lowercase letters like $v$ and the corresponding pocket 
by $Y_v$. If  $e$ is a directed edge of $T_X$ from $\del_{-}(e) = v$ to $\del_{+}(e) = w$, then we denote by $\mu_e$ the unique 
trivial product of leaves whose bottom leaf is contained in $\partial_+ Y_v$ and whose top leaf is contained in $\partial_- Y_w$. Although $\mu_e$ is a trivial product of depth $n-k$ leaves, as the notation suggests, we will often treat it as a single isolated leaf.
For simplicity, we further denote $Y_{\del_+ e}$ as $Y_+(e)$ and $Y_{\del_- e}$ as $Y_-(e)$. 
Note that, similarly to how we are now indexing subpockets of $\wt X$ using vertices of $T_X$, we have dropped the $\wt \cdot$ for ease of notation.
}
\end{remark}

For each vertex $v$ of $T_X$, the pocket $Y_v$ is the lift to $\wt M$ of a pocket of $M \ssm \mc F_{\le n-k}$ to which the inductive hypothesis applies. That is, there is a finite set of choices of universal circles  
for the restriction of $\mc F$ to the projection of $Y_v$ to $M$
such that each universal circle for a transverse flow is isomorphic to one in this finite set. For each of the finitely many vertices of $G_X$ (i.e. for each $\pi_1(X)$--orbit of vertex of $T_X$) we choose such a universal circle  $\uc_{Y_v}$, so that the collection of chosen circles is $\pi_1(X)$--equivariant. We further suppose that for our given fixed choice there is a \define{compatible} almost transverse pseudo-Anosov flow $\varphi$. This means that for each vertex $v$ of $T_X$, 
$\uc_{Y_v}$ is isomorphic to $\uc^\phi_{Y_v}$ as a universal circle. 

\smallskip

For each edge $e$ of $T_X$, the leaves in the product $\mu_e$ have uniformly bounded Hausdorff distance and are therefore quasi-isometric. Hence, we may canonically identify their ideal boundaries with a single circle that we denote by $\partial \mu_e$. Moreover, we note that for any almost transverse pseudo-Anosov flow, the shadow of $\mu_e$ is equal to the shadow of any of its leaves by \Cref{prop:products}.

\begin{remark}[Notation; indexing universal circles and monotone maps]
\emph{
Having fixed the universal circles $\uc_{Y_v}$ for each vertex $v$ of $T$, recall that there are monotone maps $\pi_{Y_\pm(e)\mu_e}$, which we denote now simply by 
\[
\pi_{\pm e} \colon \uc_{Y_\pm(e)} \to \partial \mu_e.
\]
Recall that the main structural results about these maps is \Cref{lem:gap_finiteness}, which controls the points whose preimages are nontrivial gaps.
}
\end{remark}

\subsection{Stations}
Fix $g \in \pi_1(X)$ and denote by $T^g_X$ the subtree of $T_X$ fixed pointwise by $g$. Note that for each $e \subset T^g_X$, the element $g$ fixes each leaf in $\mu_e$ and therefore acts with north-south dynamics on $\partial \mu_e$ with a unique $g$--attractor $a_g \in \partial \mu_e$. We say that $\mu_e$ has a \define{$\pm$ $g$--gap} (also called a $g$-gap from above/below) if for the monotone map
\[
\pi_{\pm e} \colon \uc_{Y_\pm(e)} \to \partial \mu_e,
\]
$\pi_{\pm e}^{-1}(a_g) \subset \uc_{Y_\pm(e)}$ is a nontrivial gap that we denote by $I_g^\pm$ (as in \Cref{cor:next}). Note that the endpoints of $I_g^\pm$ must be $g$--attractors but there may be finitely many additional $g$--attractors in its interior. 

\begin{definition}[Stations]
A non-backtracking
edge path $\gamma = e^{\epsilon_1}_1 \cdots e^{\epsilon_k}_k$ in $T^g_X$ is called a \define{$g$--station} if $\mu_{e_1}$ has a $-\epsilon_1$ $g$--gap and $\mu_{e_k}$ has an $\epsilon_k$ $g$--gap.

A \define{station} is a $g$--station for some $g \in \pi_1(X)$.
\end{definition}

Note that with this definition, if $\gamma$ is a $g$--station, then so is its reverse $\ol \gamma$. In this definition, $\epsilon_i = \pm 1$ and $e_i^{-1}$ means $e_i$ with its opposite orientation (from the natural one induced from $T_X$).
So, in particular, if $\gamma$ is a directed path in $T_X$, then each $\epsilon_i =+1$, the first $g$--gap is from below, and the last $g$--gap is from above. 



For the simplest example of a $g$--station, consider an edge path $\gamma=e$ in $T_X$ of length 1 such that the shadow $\Omega_{\mu_e}$ of $\mu_e$ has a $g$--invariant frontier chain containing a stable frontier component $l^s$ and an unstable frontier component $l^u$. The image of this frontier chain in $\partial \mu_e$ is fixed by $g$ and we suppose that it is the $g$--attractor $a_g \in \partial \mu_e$
Here, we are working in the flow space of the fixed compatible flow $\varphi$. Finally, further suppose that the shadow of $Y_+(e)$ contains $l^u$ and that the shadow of $Y_-(e)$ contains $l^s$. 

It follows that $\mu_e$ has a $g$-gap from both above and below and hence $\gamma$ is a $g$--station. In this case, the $g$--fixed point in $l^s$ corresponds to a closed orbit in the projection of $Y_-(e)$ to $M$ that is homotopic to the inverse of the closed orbit in the projection of $Y_+(e)$ to $M$ that corresponds to the $g$--fixed point in $l^u$. That is, these closed orbits are anti-homotopic. 

\smallskip
Note that if $\gamma$ is a $g$--station, then for $h \in \pi_1(X)$, $h\gamma$ is a $hgh^{-1}$--station, where the sign of the station is preserved since $\pi_1(X)$ preserves the orientation on edges of $T_X$. 

Fix a $g$--station $\gamma$. For simplicity, we will denote by $P_0(\gamma)$ its initial pocket and $P_1(\gamma)$ its terminal pocket (using the orientation determined by $\gamma$). 
We also denote by $\mu_0$ the products of depth $n-k$ leaves associated to the first edge of $\gamma$ and $\mu_1$ the products of depth $n-k$ leaves associated to the last edge of $\gamma$. We denote the monotone maps by
\[
\pi_{0 \gamma} \colon \uc_{P_0(\gamma)} \to \partial \mu_{0}
\]
and  
\[
\pi_{1 \gamma} \colon \uc_{P_1(\gamma)} \to \partial \mu_{1}.
\]

By the definition of a $g$--station, each of these 
monotone maps has a $g$--invariant gap $I_g^{0/1}$ that maps to the unique $g$--attractor $a_g^{0/1} \in \partial \mu_{0/1}$.
As above, the endpoints of each $I_g^{0/1}$ are $g$--attractors, and each contains at most finitely many other $g$--attractors. An interval in $I_g^{0/1}$ that meets $g$--attractors exactly at its endpoints will be called a \define{$g$--petal}. The reason for this terminology is that such an interval maps to a (usually not embedded) loop in the hyperbolic boundary of $\wt M$ under a natural Cannon-Thurston map, a fact we do not need---we mention it only as a potential mnemonic. 

The next result is the main place where \Cref{lem:gap_finiteness} is used.

\begin{lemma}[Station finiteness] \label{lem:station_finite}
For any $g \in \pi_1(X)$, there are finitely many $g$--stations. Moreover, up to the action of $\pi_1(X)$, there are finitely many stations. 
\end{lemma}

\begin{proof}
As before, let $\varphi$ be any almost transverse pseudo-Anosov flow that is compatible with our fixed choices of universal circles for the $Y_v$. 

Now apply \Cref{lem:gap_finiteness} to conclude that every nontrivial gap is uniquely associated to a finite set of loops. This gives a finite set of elements $g \in \pi_1(X)$ such that any station is a $hgh^{-1}$--station for some $h \in \pi_1(X)$. 

Hence, it remains to show that for any such $g\in \pi_1(X)$, there are finitely many $g$--stations. This will follow from the fact that there are finitely many leaves of $\orb^{s/u}$ fixed by $g$.
To each such $g$--station $\gamma$, using the notation before the lemma statement, the unique $g$--attractors in $\partial \mu_0$ and $\partial \mu_1$ 
determine frontier chains in $\orb_\phi$ with the property that each contains a frontier component (say $l_0$ and $l_1$, respectively) such that $l_0$ is in the shadow of the pocket $P_0(\gamma)$ and $l_1$ is in the shadow of the pocket $P_1(\gamma)$. These leaf slices are not necessarily unique, but they exist because the shadow of a pocket contains the shadow of its boundary leaves and because the gaps corresponding to the $g$--attractors are nontrivial (and therefore the corresponding frontier chains do not agree).

Note that each of $l_0$ and $l_1$ are $g$--invariant (in fact, each endpoint in $\partial l_i$ is a $g$--attractor in $\partial \orb_\phi$). There are at most finitely many such leaf slices; it therefore suffices to show that $\gamma$ is determined by the pair $l_0,l_1$ as well as a choice of side of each in $\orb$. We will do this by showing that $l_0$, plus a side choice, determines the pocket $P_0(\gamma)$ (and similarly for $l_1$). Since $T_X$ is a tree, this determines $\gamma$. 

For this, let us focus on $l_0$ and for concreteness, suppose that $l_0$ is stable. (The unstable case is similar.) Fixing a side of $l_0$, consider all the $\mu_e$ with $e \subset T_X$ such that $l_0$ is a frontier component of the shadow of $\mu_e$ and $\mu_e$ lies on the chosen side.  All these $\mu_e$ are ordered by the flow and hence they lie along a \emph{directed} path in $T_X$. If $v$ is the vertex of $T_X$ associated to $P_0(\gamma)$, we know that this directed path is actually based at $v$ and so we can find $v$ by taking the first vertex of this path. 

This shows that $l_0$ determines $v$ and hence the map from $g$--stations to pairs of leaves given by $\gamma \mapsto (l_0,l_1)$ is injective. Since the codomain is finite, so is the domain, and the proof is complete.
\end{proof}

\subsection{Plans}

We define another object that will be important in the proof of our main result.

\begin{definition}[Plans]
A \define{plan} for a $g$-station $\gamma$ is a pair of functions. One function takes as input a $g$--petal of $I_g^0$ and outputs a single $g$--attractor of $I_g^1$. The other takes a $g$--petal of $I_g^1$ and outputs a single $g$--attractor of $I_g^0$.

In general, a \define{plan} is a plan for each station that is equivariant under the natural $\pi_1(X)$--action. 
\end{definition}

A plan might be thought of as instructions for traversing the parts of $ E_\infty$ that correspond to stations. We will soon see that each plan determines one or zero universal circles.
The following corollary is immediate from the finiteness of stations and the finiteness of $g$--attractors in each circle. 

\begin{corollary}[Plan finiteness]\label{finite_plans}
There are finitely many plans.
\end{corollary}

A simplified perspective on the proof of our main finiteness theorem is that in constructing all possible universal circles for $\mc F$, ``all choices are made at stations, all of these choices are encoded by plans, and there are finitely many stations and plans."

\subsection{From plans to sections}
\label{sec:plan_to_sec}

Recall that we have fixed a pocket $X$ of depth $k+1$. Now fix a mapping torus leaf $\lambda$ of $X$ (which is nontrivial because $\mc F$ is reduced) 
and $h \in \pi_1(M)$ that uniquely stabilizes $\lambda$, as in \Cref{sec:stab}. Since $\lambda_M \subset X$, we also have that $h \in \pi_1(X)$. As before, let $x \in \partial \lambda$ be the unique $h$--attractor in $\partial \lambda$. 

Our goal in this section is to explain how a plan determines a unique section $s_x$ of the circle bundle $E_\infty|_{\wt X} \to \mc L|_{\wt X}$ with the property that $s_x(\lambda) = x \in \partial \lambda$. An outline of our strategy is the following. Note that each leaf of $\mc L|_{\wt X}$ is either 
\begin{itemize}
\item a leaf of some $Y_v$ for $v \in T^0_X$,
\item a leaf of some $\mu_e$ for an edge $e \subset T_X$, or
\item a boundary leaf of $\wt X$.
\end{itemize}
We will construct a function $s_x$ on $T^0_X$ such that $s_x(v)$ is a point of $\uc_{Y_v}$ satisfying a compatibility condition across edges of $T_X$. Namely, to each edge corresponds two pockets $Y_+(e)$ and $Y_-(e)$, and a single circle $\del \mu_e$; we require that the values of $s_x$ on $Y_\pm (e)$ project to the same point in $\del \mu_e$. In symbols,
\begin{align}\label{disc_sec}
\pi_{+e}(s_x(Y_+(e))) = \pi_{-e}(s_x(Y_-(e))).
\end{align}

To determine the value of our section at a leaf $\lambda$ in a pocket $Y_v$, we will simply project $s_x(Y_v)$ to $\del \lambda$ using the appropriate monotone map. This will determine a continuous section of $E_\infty|_{\wt X} \to \mc L|_{\wt X}$, the value of the section at leaves of $\partial_\pm \wt X$ being determined by continuity (as in \Cref{rmk:frontier_leaves}).

\begin{construction} \label{plan_to_sec}
For $x$ as above, a plan determines a well-defined section $s_x$ with $s_x(\lambda) = x$. 
\end{construction}

The construction is as follows.
Let $v_0 \in T_X$ denote the vertex such that $\lambda \subset Y_{v_0}$. The section $s_x$ is defined inductively for each $v \in T_X$, where the induction is on the distance to $v_0$ in $T_X$. For the base case, we claim that there is a unique point in $\pi^{-1}_{Y_{v_0}\lambda}(x) \subset \uc_{Y_{v_0}}$, which we then set to be $s_x(v_0)$. For the argument, 
recall that we have fixed a compatible, almost transverse pseudo-Anosov flow $\varphi$ so that, in particular, $\uc_{Y_{v_0}}$ is isomorphic to $\uc^\varphi_{Y_{v_0}}$. 
By \Cref{prop:north_south}, the element $h$ fixed at the start of this subsection
has a unique attractor $\wh x$ in $\partial \orb_\phi$ and so $(\pi^\phi_{\lambda})^{-1}(x) = \{\wh x\}$. 
Hence, $\pi^{-1}_{Y_{v_0}\lambda}(x) = \{\pi^\varphi_{Y_{v_0}}(\wh x)\}$ and the claim follows. 

Note a subtle point: despite the fact that above we used the flow $\phi$ to argue that $s_x(v_0)$ is well-defined, its definition (as the lone point in $\pi^{-1}_{Y_{v_0}\lambda}(x)$) is flow-independent.

Now suppose $s_x$ is defined for all vertices of $T_X$ whose distance to $v_0$ is less than $d$. Let $v$ be a vertex at distance $d$ and let $e$ be the terminal edge of a length $d$ path from $v_0$ to $v$. For simplicity, assume that the orientation of $e$ coming from this path agrees with the natural orientation (the other case is similar).  Hence, $v = \partial_+ e$, and we set $w = \partial_- e$. Since $s_x(w)$ is already defined, we can set $s_x(\mu_e) = \pi_{Y_w\mu_e}(s_x(w)) \in \partial \mu_e$. 
\smallskip

If $\pi_{Y_v \mu_e}^{-1}(s_x(\mu_e))$ is a single point, then we set this point to be $s_x(Y_v)$. In this case, we say that the plan \emph{was not consulted}.

Otherwise, $\pi_{Y_v \mu_e}^{-1}(s_x(\mu_e))$ is a nontrivial gap.  From \Cref{lem:gap_finiteness} and the remark following it, this implies that $s_x(\mu_e) \in \partial \mu_e$ is a $g$--attractor for a unique primitive $g$ stabilizing $\mu_e$. 

\smallskip
A \define{relevant $g$--station $\gamma$ for the pair $(v_0,v)$} is a $g$--station $\gamma$ contained in the path in $T_X$ from $v_0$ to $v$ and ending at $v$ (i.e. $P_1(\gamma)= v$), so $s_x(P_0(\gamma)) \in I_g^0$. That is, the value of $s_x$ at the start of the station is immediately collapsed to the unique $g$--attractor through the next monotone map.

Now if there is no relevant $g$-station contained along the path from $v_0$ to $v$,
then we define $s_x(Y_v)$ to be the leftmost point of $\pi_{Y_v \mu_e}^{-1}(s_x(\mu_e))$, which is necessarily a $g$--attractor as recalled above.

So assume that $\gamma$ is the
largest relevant $g$--station contained in the path from  $v_0$ to $v$ in $T_X$. By definition, $s_x(P_0(\gamma)) \in I_g^0$ and so $s_x(P_0(\gamma))$ is contained in a petal in $I_g^0$ (if it happens to be a $g$--attractor, then it is in two petals and we take the leftmost). The fixed plan then associates to this petal a unique $g$--attractor in $I_g^1 \subset \uc_{Y_v}$. We take this point for $s_x(Y_v)$. 

This defines $s_x$ for all vertices of distance $d$ from $v_0$. This defines $s_x$ uniquely from a plan.

\subsection{From flows to plans}
In this short subsection, we show that any transverse almost pseudo-Anosov flow determines a plan, which we call its associated \define{flow plan}. In the next section, we will use the resulting plan to recover the universal circle itself. 

\begin{construction} \label{flow_to_plan}
A flow $\phi$ determines a well-defined plan.
\end{construction}

Fix a $g$--station $\gamma$. We must use $\varphi$ to define a function from the petals of $I_g^0 \subset \uc_{P_0(\gamma)}$ to the $g$--attractors of $I_g^1 \subset \uc_{P_1(\gamma)}$. 
Recall that $P_0(\gamma)$ and $P_1(\gamma)$ are the pockets at the beginning and end of the $\gamma$, respectively, and $I_g^{0/1}$ are the gaps associated to the unique $g$--attractors. 

Let $p \subset I_g^0$ be a petal. We call the open interval
\[
U ={(\pi^{\phi}_{P_0(\gamma)})}^{-1}(\mathrm{int}(p)) \subset \partial \orb_\phi
\]
a \define{pre-petal} in $\partial \orb_\phi$. 

If there is a frontier component $l_1$ of $\Omega_{P_1(\gamma)}$ such that 
\begin{itemize}
\item $l_1$ separates $\Omega_{P_1(\gamma)}$ from $U$, and 
\item $\pi^\phi_{\mu_1}(l_1)$ is the unique $g$--attractor in $\partial {\mu_1}$,
\end{itemize}
then we set the value of the plan for $\gamma$ at the petal $p$ to be $\pi^\phi_{P_1(\gamma)}(l_1)$. Here, recall that $\mu_1$ is the product of depth $n-k$ leaves corresponding to the last edge of $\gamma$. Hence, $\mu_1$ contains a $g$--invariant boundary leaf of $P_1(\gamma)$. Note that by the second item, $l_1$ is $g$--invariant and its endpoints map to the same $g$--attractor under $\pi^\phi_{P_1(\gamma)}$. Hence, in this case, we satisfy the definition of a plan.

Otherwise, we take the value of the plan for $\gamma$ at the petal $p$ to be the leftmost $g$--attractor of $I_g^1 \subset \uc_{P_1(\gamma)}$. This last choice is arbitrary, but as we will see below, it doesn't come into play in what follows.

\subsection{Flow plans give flow space sections}
\label{sec:flow_plan}
Recall that $x \in \partial \lambda$ is fixed as in \Cref{plan_to_sec} as the attractor of $h$, where $\lambda$ is uniquely stabilized by $h$. By \Cref{prop:north_south}, for any transverse flow $\varphi$, $h$ has a unique attractor $\wh x$ in $\partial \orb_\phi$ and so $\pi^\phi_\lambda(\wh x) = x$. 

By \Cref{flow_to_plan}, $\phi$ determines a plan, and by \Cref{plan_to_sec}, this plan determines a section $s_x$. We call $s_x$ the \define{section determined by the flow plan for $\phi$.}

\begin{lemma} \label{lem:right_plan}
For any transverse almost pseudo-Anosov flow $\phi$, let $s_x$ be the section determined by the flow plan for $\phi$. Then for any leaf $\mu$ of $\wt X$:
 \[
 s_x(\mu) = \pi^\phi_\mu(\wh x).
 \] 
\end{lemma}

By the discussion surounding \Cref{disc_sec}
(and \Cref{rmk:frontier_leaves}), it suffices to show that for each $v \in T^0_X$,  $s_x(Y_v) = \pi^\phi_{Y_v}(\wh x)$.
Recall that $s_x$ was defined inductively with $s_x(Y_{v_0}) = x$, where the induction was over the distance from a vertex $v$ to the base vertex $v_0$, where $Y_{v_0}$ contains $\lambda$. Hence, the base case holds by construction; so suppose it holds for all vertices whose distance from $v_0$ is less than $d$. 

Also, recall by \Cref{lem:good_to_move}, the shadow of every leaf not parallel to $\lambda$
has a boundary leaf that separates it from $\wh x$ in $\ol \orb_\phi$. Hence, the same is also true for the shadow of every $Y_v$, other than $Y_{v_0}$. We will use this fact to apply the Moving Away Lemma \ref{lem:moving_saturated}.

So, as in the proof of \Cref{plan_to_sec}, let $v$ be a vertex whose distance to $v_0$ is $d$. We again let $e$ be the last edge of the length $d$ path from $v_0$ to $v$ and again assume for simplicity that the orientations on $e$ determined by this path and $T_X$ agree. Hence, $v = \partial_+ e$ and we let $w = \partial_- e$. Note that by induction $s_x(w) = \pi^\phi_{Y_w}(\wh x)$ and so in particular, $s_x(\mu_e) = \pi^\varphi_{\mu_e}(\wh x)$. 

\smallskip
Recall that if $\pi_{Y_v \mu_e}^{-1}(s_x(\mu_e))$ is a single point, then the plan was not consulted and this point was declared to be $s_x(Y_v)$. Since  $s_x(\mu_e) = \pi^{\varphi}_{\mu_e}(\wh x)$, this point must be $\pi^\varphi_{Y_v}(\wh x)$ as required. 

\begin{remark}[No perfect fits] \label{rmk:npfs}
\emph{
Since $\mu_e$ separates $\lambda$ from $Y_v$,
we are in the setting of \Cref{cor:next}. Hence, if $\mc F$ is almost transverse to 
\emph{some} pseudo-Anosov flow 
without perfect fits, then the `moreover' claim of \Cref{cor:next} states that 
$\pi_{Y_v \mu_e}^{-1}(s_x(\mu_e)) = \pi_{Y_v \mu_e}^{-1}(\pi^{\phi}(\wh x))$ is a single point. Hence, in this case, the plan is never
consulted. 
}

\end{remark}

Otherwise, $\pi_{Y_v \mu_e}^{-1}(s_x(\mu_e)) = \pi_{Y_v \mu_e}^{-1}(\pi^{\phi}(\wh x))$ is a nontrivial gap.  From \Cref{cor:next},
this implies that $s_x(\mu_e) \in \partial \mu_e$ is a $g$--attractor for a unique (primitive) $g$ stabilizing $\mu_e$.

The next claim holds because of the special nature of the section $s_x$:
\begin{claim}[Relevant $g$--station] \label{relevant}
There is a unique largest relevant $g$--station $\gamma$ for the pair $(v_0,v)$.
 \end{claim}

To simplify notation, set $\mu = \mu_e$.
 
The main point is that a relevant $g$-station exists, so we will produce one. Consider the frontier component $l_\mu$ of $\Omega_\mu$ that separates it from $\wh x$. Working backwards along the unique embedded path from $v_0$ to $v$, there must have been some first $w \in T^0_X$ so that
 the shadow of $Y_w$ contains $l_\mu$ and the next edge $f$ (along the path from $w$ to $v$) corresponds to a trivial product of leaves $\mu_f$ meeting the boundary of $Y_w$ 
 whose shadow also has $l_\mu$ as a frontier component.
Note that either the limit set of the shadow of $Y_w$ contains $\wh x$ (which occurs if and only if  $w = v_0$) or the frontier component that separates it from $x$ also separates $\Omega_\mu$ from $x$ by \Cref{lem:moving_saturated}.
 
This shows that the path from $w$ to $v$ is a relevant $g$--station for the pair $(v_0,v)$. Here, we are using the inductive hypothesis that $s_x(Y_w) = \pi^\phi_{Y_w}(\wh x)$ to establish that the station is relevant. 
The claim now follows by taking the maximal relevant $g$-station. We will denote this station by $\gamma$.

\smallskip

With $\gamma$ now fixed, let $u = P_0(\gamma)$ be its initial vertex and $b$ its initial edge; by construction $v = P_1(\gamma)$.
Note that in the path from $v_0$ to $v$, $u$ comes before $w$, and so again by \Cref{lem:moving_saturated} 
the frontier component $l_{b}$ of the shadow of $\mu_{b}$ that separates it from $\wh x$ also separates the shadow of $\mu_f$ (and hence the frontier component $l_\mu$) from $\wh x$.

The following claim is important for the plan to build the correct section:

\begin{claim}
The point $s_x(Y_{u})$ is \emph{not} a $g$--attractor of the gap $\pi_{Y_{u}\mu_{b}}^{-1}(s_x(\mu_{b})) \subset \uc_{Y_{u}}$. As a consequence, the frontier component $l_b$ of the shadow of $\mu_{b}$ that separates it from $\wh x$ lies in the shadow of $Y_u$.
\end{claim}

This implies $s_x(Y_{u})$ is contained in the interior of a unique petal of the gap $\pi_{Y_{u}\mu_{b}}^{-1}(s_x(\mu_{b})) \subset \uc_{Y_{u}}$.

We will again use that by our inductive hypothesis: $s_x(Y_u) = \pi^\phi_{Y_u}(\wh x)$.

\begin{proof}
First note that the second statement follows from the first. This is because if $l_b$ were not in $\Omega_{Y_u}$, then it would be in its frontier and it would be the unique frontier component separating it from $\wh x$. But $l_b$ is $g$--invariant and its endpoints are $g$--attractors, implying that $s_x(Y_u) = \pi^\phi_{Y_u}(\wh x) = \pi^\phi_{Y_u}(\partial l_b)$ is a $g$--attractor, contradicting the first statement.

For the first statement, suppose otherwise. Since $\wh x$ is not the endpoint of a leaf, there must be a frontier component $l_u$ of $\Omega_{Y_u}$ that separates it from $\wh x$, such that the endpoints of $l_u$ are $g$-attractors.. 
(Note that if $u=v_0$, then $s_x(Y_u) =x$ is not a $g$--attractor.)
Just as in the proof of \Cref{relevant}, we can move backwards along the unique embedded path in $T_X$ from $v_0$ to $u$ to find a vertex $u'$ such that the shadow of $Y_{u'}$ contains $l_u$. Just as in that setting, this shows that the path from $u'$ to $v$ is a relevant $g$--station, contradicting that $\gamma$ is the largest such.
\end{proof}

\begin{claim}
Let $p$ be the unique petal of $\uc_{Y_u}$ whose interior contains $s_x(Y_u) = \pi^\phi_{Y_u}(\wh x)$. The frontier component $l_b$ that separates the shadow of $\mu_b$ from $\wh x$ also separates it from the pre-petal $U = {(\pi^{\phi}_{Y_u})}^{-1}(\mathrm{int}(p))$. 
\end{claim}

\begin{proof}
Note that $\wh x$ is contained in the pre-petal $U$.

By the previous claim, the shadow of $Y_u$ contains $l_b$. So the endpoints of $l_b$ project to distinct $g$--attractors in $\uc_{Y_u}$. On the other hand, a pre-petal does not contain any point in the boundary of the flow space that projects to a $g$--attractor in $\uc_{Y_u}$. Hence $U$ must lie entirely in the component of $\del \orb \ssm\{\text{endpoint of $l_b$}\}$ that contains $\wh x$. This implies the claim.
\end{proof}

We are now ready to conclude our proof by establishing that $s_x(Y_v) = \pi^\phi_{Y_v}(\wh x)$. Indeed, by the previous claim and the Moving Away Lemma (\Cref{lem:moving_away}), $l_b$ separates the shadow of $\mu$ from both $\wh x$ and the pre-petal $U$. Hence, the same is true of the frontier component $l_\mu$ of $\Omega_\mu$ that separates it from $\wh x$. By \Cref{lem:moving_saturated}, $l_\mu$ is also a frontier component of $\Omega_{Y_v}$. 

Hence, from the definition of the flow plan given after the statement of \Cref{flow_to_plan}, we conclude that 
\[
s_x(Y_v) = \pi^\phi_{Y_v}(l_\mu) = \pi^\phi_{Y_v}(\wh x),
\]
and the proof is complete. 

\medskip
Finally:

\begin{proof}[Proof of \Cref{prop:first_col}]
Throughout the proof we fix the lifted mapping torus leaf $\lambda$ and point $x\in \del\lambda$, which is the unique attractor for a group element uniquely stabilizing $\lambda$, as we have throughout this section.

Each of the finitely many plans (\Cref{finite_plans}) determines a unique section (\Cref{plan_to_sec}), which by \Cref{lem:good_order} is contained in at most one universal circle of $\mc F|_X$ with a minimal action. Let $\uc_1, \ldots, \uc_l$ be this finite collection of universal circles for $\mc F|_X$.

Now if $\varphi$ is a pseudo-Anosov flow on $M$ that is almost transverse to $\mc F$, then it determines a plan via \Cref{flow_to_plan}. This plan determines a section 
\[
s_x\colon \mc L\to E_\infty
\]
via \Cref{plan_to_sec}.
On the other hand, $\phi$ has a flow space universal circle $\uc^\phi_X$ which has a section defined by 
\[
\mu \mapsto \pi_\mu^\phi(\wh x), 
\]
where $\wh x$ is the unique point of $\partial \orb_\phi$ such that $\pi^\phi_\lambda(\wh x) =x \in \partial \lambda$. By \Cref{lem:right_plan}, these two sections agree. Since $\uc_X^\phi$ is minimal by \Cref{prop:minimal}, it must be isomorphic to one of the $\uc_i$ determined above. 


Finally, if $\mc F$ is almost transverse to a pseudo-Anosov flow without perfect fits, 
then it suffices to show that any plan 
produces the same section $s_x$, and so there is a unique universal circle $\uc_1$ in our list, up to isomorphism. This follows from the fact that the moreover statement in \Cref{cor:next} implies that when such a transverse flow exists, the plan is never consulted when the section is built in \Cref{plan_to_sec} (see \Cref{rmk:npfs}).
In other words, $s_x$ is unique, as is the universal circle $\uc_1$ it determines. 
Hence $\mc F|_X$ has a unique flow space universal circle, as required.
\end{proof}

\subsection{Proofs of \Cref{th:main_tft} and \Cref{th:main_npfs}}

Apply \Cref{prop:first_col} to the case where $k=n+1$. This produces a finite set of universal circles that contains the isomorphism class of the flow space boundary of any almost transverse pseudo-Anosov flow. Since two flows with the same flow space boundary (considered as a circle action) are orbit equivalent when the underlying manifold is atoroidal by \cite[Theorem 1.5]{BFM}, this completes the proof of the transverse finiteness theorem (\Cref{th:main_tft}).

If there is an almost transverse pseudo-Anosov flow without perfect fits, then there is at most one universal circle on our list. This implies that there is at most one almost transverse pseudo-Anosov flow. This completes the proof of our rigidity theorem (\Cref{th:main_npfs}).

\section{Universal circle rigidity}
\label{sec:uc_npfs}

In this section, we prove the following restatement of \Cref{th:main_ucs}.
\begin{theorem}
\label{th:UC_rig}
Suppose that $\mc F$ is a finite depth foliation that is almost transverse to a pseudo-Anosov flow $\varphi$ without perfect fits. Then for any minimal universal circle $\uc$ for $\mc F$, $\uc$ is isomorphic to $\partial \orb_\varphi$.
%
\end{theorem}

For simplicity, we write $\uc_1 = \uc_2$ to indicate that $\uc_1$ and $\uc_2$ are isomorphic. This notation is justified by the fact that when an isomorphism exists, it is unique.
See the discussion following \Cref{def:uc}. 
\smallskip

Toward the proof, fix such a universal circle $\uc$ of $\mc F$. Since we will follow the previous induction scheme, we first make the observation that $\uc$ induces a universal circle on each open $\mc F$--saturated set (in the sense of \Cref{rmk:uc_sat}) by following a procedure nearly identical to the one in \Cref{sec:flow_ucs}.
More precisely, 
if $X$ is the completion of an $\mc F$--saturated open set of $M$ (e.g. a pocket) with lift $\wt X$, then we define the circle
\[
\uc_{\wt X} = \ol{\bigcup_{\mu \subset \wt X} \mathrm{core}(\pi^{\uc}_\mu)} / \sim,
\]
where $\pi^{\uc}_\mu \colon \uc \to \partial \mu$ is the monotone map to the boundary of a leaf $\mu \subset \wt X$ and $\sim$ collapses the closures of the complementary intervals of $\ol{\bigcup_{\mu \subset \wt X} \mathrm{core}(\pi^{\uc}_\mu)}  \subset \uc$. Clearly, there is an induced action $\pi_1(X) \curvearrowright \uc_{\wt X}$.

As in \Cref{sec:flow_ucs}, there are also induced monotone maps $\pi^{\uc}_{\wt X} \colon \uc \to \uc_{\wt X}$ and $\pi^{\uc}_{\wt X \mu} \colon \uc_{\wt X} \to \partial \mu$ for each $\mu \subset \wt X$ and these satisfy: 
\begin{align} \label{eq:comp}
\pi^{\uc}_\mu = \pi^{\uc}_{\wt X\mu}\circ \pi^{\uc}_{\wt X}.
\end{align} 
The same holds if $\partial \mu$ is replaced by $\uc_{\wt Y}$, where $\wt Y \subset \wt X$ is the lift of a subpocket of $X$.

\smallskip
The following lemma is an easier version of \Cref{prop:minimal}. It follows from the first two paragraphs of that proof, together with the standardness of leaves of $\mc F$ (\Cref{prop:standard}) and the definition of $\uc_{\wt X}$. 

\begin{lemma}\label{lem:minimal2}
If $\mc F$ is finite depth and $X$ is the completion of an $\mc F$--saturated open set, then the action $\pi_1(X) \curvearrowright \uc_{\wt X}$ is minimal.
\end{lemma}

In fact, the same argument from the proof of \Cref{prop:minimal} implies that if $\del\orb$ is the flow space universal circle for an almost transverse pseudo-Anosov flow $\varphi$, then 
\[
(\del\orb)_{\wt X} = \uc^\varphi_{\wt X}
\]
where $\uc^\varphi_{\wt X}$ is as defined in \Cref{sec:flow_ucs}. However, we will not need this fact directly.
%
%

\subsection{Product compatible universal circles}
For this subsection, let $\mc F$ be an arbitrary taut foliation. If $X$ is an $\mc F$--saturated subset of $M$ that is a topological product (i.e. $X \cong \lambda_M \times [0,1]$), then any lift $\wt X$ of $X$ to $\wt M$ has the property that 
any two leaves $\lambda_1, \lambda_2 \subset \wt X$
have bounded Hausdorff distance (since $M$ is compact the fibers $\{\mathrm{pt.}\}\times [0,1]$ have bounded length). Hence, the nearest point projection between these leaves determines a quasi-isometry from $\lambda_1$ to $\lambda_2$
that induces a homeomorphism from $\partial \lambda_1$ to $\partial \lambda_2$.
We call this the \define{natural identification} of 
$\partial \lambda_1$ and $\partial \lambda_2$. When 
$\lambda_1 = \partial_-{\wt X}$ and $\lambda_2 =\partial_+ \wt X$, the natural 
identification is equivariant with respect to the $\pi_1(X)$--action.

A universal circle $\uc$ for $\mc F$ is \define{product compatible} if whenever $X$ is an $\mc F$--saturated topological product and $\lambda_1, \lambda_2$ are leaves of a lift $\wt X$, the monotone maps $\pi^{\uc}_{\lambda_1}$ and $\pi^{\uc}_{\lambda_2}$ are equal after applying the natural identification of $\partial \lambda_1$ and $\partial \lambda_2$. For example, \Cref{prop:products} implies that for any pseudo-Anosov flow $\varphi$ almost transverse to $\mc F$, the universal circle $\partial \orb_\varphi$ is product compatible. 

\begin{lemma}\label{lem:compatible}
A universal circle $\uc$ of a finite depth foliation $\mc F$ is product compatible.
\end{lemma}

\begin{proof}

Let $X$ be an $\mc F$--saturated topological product, let $\lambda_M = \partial_-X$, and fix an identification $X = \lambda_M \times [0,1]$ so that the leaves of $X$ are transverse to the $[0,1]$--fibers. This can be done, for example, by fixing an almost transverse pseudo-Anosov flow as in \Cref{prop:products}. For any leaf $\mu_M$ of $X$, the projection along the $[0,1]$--fibers is a covering map to $\lambda_M$ that identifies $\pi_1(\mu_M)$ with a subgroup of $\pi_1(\lambda_M)$. If we fix a lift $\wt X$ of $X$ to $\wt M$ and use $\lambda$ and $\mu$ to denote the corresponding lifts of $\lambda_M$ and $\mu_M$, then the natural identification $\partial \mu= \partial \lambda$ identifies the boundary action $\pi_1(\mu_{M}) \curvearrowright \partial \mu$ with the action $\pi_1(\mu_M) \curvearrowright \partial \lambda$ induced by the subgroup $\pi_1(\mu_M) \le \pi_1(\lambda_M)$. 

Hence, it suffices to show that there is at most one $\pi_1(\mu_M)$--equivariant monotone map $\uc \to \partial \lambda$. This follows quickly from the fact that the action $\pi_1(\mu_{M}) \curvearrowright \partial \mu = \partial \lambda$ is minimal (\Cref{prop:standard}). Indeed, because of minimality, the core $C \subset \uc$ of any such monotone map $\uc \to \partial \lambda$ is uniquely determined as the smallest closed $\pi_1(\mu_M)$--invariant subset of $\uc$. Moreover, the map on the core $C$ is determined by equivariance:
for each nontrivial $g\in \pi_1(\mu_M)$, $g$ has a unique attracting point $a_g \in \partial \lambda$ and hence at most two (local) attracting points in the core $C$, which must get mapped to $a_g$. Since such attracting points are dense in $\partial\lambda$ and a monotone map is determined by the restriction to its core, the lemma follows. \qedhere


%

\end{proof}

\subsection{Proof of \Cref{th:UC_rig}}
With this setup, we can start the proof.

\begin{proof}[Proof of \Cref{th:UC_rig}]
We first assume that the finite depth foliation $\mc F$ is reduced.

We will show that for any pocket $X$ of $M\ssm\mc F_{\le n-k}$, there is an isomorphism $\uc_{\wt X} = \uc^\phi_{\wt X}$. Applying this to the case $X = M$ proves the theorem, since by definition $\uc^\phi_{\wt M} = \partial \orb_\phi$. The proof is by induction on $k$ and has the exact same setup as in \Cref{sec:setup}, which the reader can consult as a reminder. Recall that by \Cref{lem:compatible}, $\uc$ is product compatible.

Within the inductive step, there will be another induction. To avoid confusion, we will refer to the induction on $k$ as the \emph{depth induction}.

In the base case of the depth induction, $X$ has top depth and is therefore a mapping torus component. Any two mapping torus leaves of $X$ cobound a product, and $\uc$ is product compatible. Hence for any lift $\wt X$, and any mapping torus leaves $\lambda_1$ and $\lambda_2$ in $\wt X$, we have that $\mathrm{core}(\pi^{\uc}_{\lambda_1}) =  \mathrm{core}(\pi^{\uc}_{\lambda_2})$. We conclude that for any mapping torus leaf $\lambda$ of $\wt X$:
\[
\uc_{\wt X} = \mathrm{core}(\pi^{\uc}_{\lambda})/ \sim \; = \partial \lambda = \uc^\varphi_{\wt X},
\]
where in the first equality
we have used the fact that if $\mu$ is a boundary leaf of $\wt X$, then $\core(\pi^{\uc}_{\mu})\subset \core(\pi^{\uc}_{\lambda})$. 

Now assume that $X$ is a pocket with lift $\wt X$. For the inductive hypothesis of the depth induction, we assume $\uc_{\wt Y}  = \uc^\phi_{\wt Y}$ for any subpocket $\wt Y \subset \wt X$ of higher depth. 

Just as in \Cref{sec:plan_to_sec}, fix a mapping torus leaf $\lambda$ of $\wt X$ and a group element $h \in \pi_1(X)$ that uniquely stabilizes $\lambda$. As at the start of \Cref{sec:flow_plan}, if $x$ is the attracting fixed point for the action of $h$ on $\partial \lambda$, there is a unique $\wh x \in \partial \orb_\varphi$ so that $\pi_\lambda^\varphi(\wh x) = x$. 


Let $s_x$ be the section of $E_\infty|_{\wt X} \to \LL_{\wt X}$ determined from the flow plan for $\varphi$. By \Cref{lem:right_plan}, we have
\[
s_x(\mu) = \pi^\varphi_\mu(\wh x)
\]
for any leaf $\mu$ of $\wt X$.
On the other hand, the main idea of what follows will be that, as recorded in \Cref{rmk:npfs}, $s_x$ is a canonical, flow-independent object. 

We can also construct a section for $\uc$ associated to $x$, although it is not obviously canonical: let $y\in \uc$ be any point with $\pi^{\uc}_\lambda(y) = x \in \partial \lambda$, and let $c_x$ be the section of $E_\infty|_{\wt X} \to \LL_{\wt X}$ determined by $y$, that is
\[
c_x(\mu) = \pi^{\uc}_\mu (y). 
\]
We will show that these two sections are equal: $s_x(\mu) = c_x(\mu)$ for all leaves $\mu$ of $\wt X$. Since both $\uc_{\wt X}$ and $\uc^\varphi_{\wt X}$ have minimal $\pi_1(X)$--actions (by \Cref{lem:minimal2} and \Cref{prop:minimal} respectively), 
we will conclude that these two universal circles are isomorphic by \Cref{lem:good_order}. This will complete the proof in the case when $\mc F$ is reduced.

\smallskip
Toward this goal, recall the dual tree $T_X$ defined in \Cref{sec:setup} as well as the notation introduced there. As in the discussion surrounding \Cref{disc_sec} (applied just as in the proof of \Cref{lem:right_plan}), the depth inductive hypothesis implies that it suffices to show the following:

\begin{claim}
After identifying $\uc_{Y_v} = \uc^\varphi_{Y_v}$ using the depth inductive hypothesis, we have
\[
\pi^{\uc}_{Y_v}(y) = \pi^\varphi_{Y_v}(\wh x)
\]
for each vertex $v$ of $T_X$.
\end{claim}

Let $v_0$ be the vertex of $T_X$ such that $Y_{v_0}$ contains the leaf $\lambda$.
Similarly to the proof of \Cref{lem:right_plan}, we will induct on the distance from $v$ to $v_0$. To distinguish this from our ongoing depth induction, we refer to this as our \emph{distance induction}. 

The distance induction base case is when $v=v_0$. This case will follow from the depth induction hypothesis, together with the fact that $Y_{v_0}$ contains $\lambda$. Indeed, 
since $\pi^\varphi_\lambda$ factors through the surjection $\pi^\varphi_{Y_{v_0}}$, the uniqueness of $\wh x$ implies that
 $(\pi^\phi_{Y_{v_0}\lambda})^{-1}(x)$ is the single point $\pi^\varphi_{Y_{v_0}}(\wh x)$. 
 On the other hand, since $\uc_{Y_{v_0}} = \uc^\varphi_{Y_{v_0}}$ by the depth induction hypothesis, we know
$\pi^\phi_{Y_{v_0}\lambda} = \pi^\uc_{Y_{v_0}\lambda}$. Hence
\[
\pi^\uc_{Y_{v_0}}(y)\in(\pi^\uc_{Y_{v_0}\lambda})^{-1}(x)=(\pi^\phi_{Y_{v_0\lambda}})^{-1}(x)=\{\pi^\phi_{Y_{v_0}}(\wh x)\}.
\]
The inclusion above follows from \cref{eq:comp}. This establishes the base case of the distance induction.

Now suppose that the distance from $v_0$ to $v$ is $d\ge 1$, and that the claim holds for vertices at distance $d-1$. Consider the nonbacktracking path from $v_0$ to $v$ in $T_X$. Let $w$ be the last vertex along this path before $v$ and let $e$ be the edge from $w$ to $v$. Since the distance from $v_0$ to $w$ is $d-1$, the distance induction hypothesis applied at $Y_w$ implies that $\pi^{\uc}_{Y_w}(y) = \pi^\varphi_{Y_w}(\wh x)$. Since $\uc_{Y_w} = \uc^\varphi_{Y_w}$, these project to a single point $z \in \partial \mu_e$ under the common monotone map $\pi^\uc_{Y_w\mu_e} = \pi^\varphi_{Y_w\mu_e}$. Then
\[
\pi^{\uc}_{Y_v}(y) \in (\pi^\uc_{Y_v\mu_e})^{-1}(z) = (\pi^\phi_{Y_v\mu_e})^{-1}(z) = \{\pi^\varphi_{Y_y}(\wh x) \}.
\]
Here the inclusion follows from \cref{eq:comp}, the first equality follows from the depth induction hypothesis (i.e. $\uc_{Y_v} = \uc^\varphi_{Y_v}$), and the second equality is \Cref{cor:next}. This completes the distance induction, establishes the claim, and completes the proof when $\mc F$ is reduced.

\medskip

It remains to consider the case where $\mc F$ is not reduced. In this case, we apply \Cref{lem:reduction} to produce a reduced finite depth foliation $\mc F_r$ obtained from $\mc F$ by refoliating product pockets and thickening isolated leaves. Recall that $\uc$ is a product compatible, minimal universal circle for $\mc F$ and $\varphi$ is the pseudo-Anosov flow without perfect fits that is almost transverse to $\mc F$. By \Cref{lem:reduction}, $\varphi$ is also almost transverse to $\mc F_r$, and it is easy to see that $\uc$ induces a product compatible, minimal universal circle $\uc_r$ for $\mc F_r$. Hence, $\uc_r = \partial \orb_\varphi$ as universal circles for $\mc F_r$. Since $\uc$ is product compatible, we can reconstruct $\uc$ from $\uc_r$ by taking the values of the monotone maps along lifted product pockets to be related by the natural identifications. Also, the monotone map of an isolated leaf can be recovered from either boundary leaf of the thickened product.
This implies that $\uc = \partial \orb_\varphi$ as universal circles of $\mc F$ and completes the proof. 
\end{proof}

\bibliographystyle{alphaurl}
\bibliography{bibliography}

\end{document}